\documentclass{article}
\usepackage[preprint,nonatbib]{neurips_2020}
\usepackage{mathtools}
\DeclarePairedDelimiter\abs{\lvert}{\rvert}

\usepackage[utf8]{inputenc} % allow utf-8 input
\usepackage[T1]{fontenc}    % use 8-bit T1 fonts
\usepackage{hyperref}       % hyperlinks
\usepackage{url}            % simple URL typesetting
\usepackage{booktabs}       % professional-quality tables
\usepackage{amsfonts}       % blackboard math symbols
\usepackage{nicefrac}       % compact symbols for 1/2, etc.
\usepackage{microtype}      % microtypography
\usepackage{adjustbox}

\usepackage{graphicx}
\usepackage{subfigure}
\usepackage[font=small,labelfont=bf]{caption}

\usepackage{wrapfig}

\newcommand{\E}{\mathbb{E}}

\usepackage[table,xcdraw]{xcolor}
\usepackage{lineno}
\usepackage{color}
\usepackage{comment}
\usepackage{algorithm}
\usepackage{algorithmic}
\usepackage{amsthm}

\usepackage[nocompress]{cite}

\usepackage{flushend}
\DeclareMathOperator*{\argmin}{argmin}

\usepackage{balance}
\usepackage{enumitem}

\usepackage{capt-of}

\newcommand{\eat}[1]{}
\newcommand{\techrep}[1]{}

\newtheorem{theorem}{Theorem} 
\newtheorem{proposition}[theorem]{Proposition}

\newtheorem{definition}{Definition}

\usepackage{array}
\usepackage{xspace}
\usepackage{makecell}

\newcommand{\stitle}[1]{\noindent{\bf #1}}

\newcommand\new[1]{\textcolor{black}{#1}}

\newcommand\attrib[1]{{\tt \small #1}\xspace}
\begin{document}
% \let\WriteBookmarks\relax
% \def\floatpagepagefraction{1}
% \def\textpagefraction{.001}

% \title [mode = title]{Fair and Diverse Allocation of Scarce Resources}                      
% \tnotemark[1,2]

\title{Fair and Diverse Allocation of Scarce Resources}
% \author{%
%   Hadis Anahideh}
% \cormark[1]
% \fnmark[1]
% \ead{hadis@uic.edu}
% \address[1]{University of Illinois at Chicago}

% \author{%
%   Lulu Kang}
% \cormark[1]
% \fnmark[1]
% \ead{lkang2@iit.edu}
% \address[2]{Illinois Institute of Technology}

% \author{%
%   Nazanin Nezami}
% \cormark[1]
% \fnmark[1]
% \ead{nnezam2@uic.edu}
% \address[1]{University of Illinois at Chicago}

\author{%
  Hadis Anahideh\\
  University of Illinois at Chicago\\
  \texttt{hadis@uic.edu} \\
  \And
  Lulu Kang\\
  Illinois Institute of Technology\\
  \texttt{lkang2@iit.edu} \\
  \And
  Nazanin Nezami\\
  University of Illinois at Chicago\\
  \texttt{nnezam2@uic.edu} \\
}

\maketitle
\begin{abstract}
% This template helps you to create a properly formatted \LaTeX\ manuscript.
% \noindent\texttt{\textbackslash begin{abstract}} \dots 
% \texttt{\textbackslash end{abstract}} and
% \verb+\begin{keyword}+ \verb+...+ \verb+\end{keyword}+ 
We aim to design a \emph{fairness-aware} allocation approach to maximize the geographical diversity and avoid \emph{unfairness} in the sense of demographic disparity. 
During the development of this work, the COVID-19 pandemic is still spreading in the U.S. and other parts of the world on large scale. 
Many poor communities and minority groups are much more vulnerable than the rest. 
To provide sufficient vaccine and medical resources to all residents and effectively stop the further spreading of the pandemic, the average medical resources per capita of a community should be independent of the community's demographic features but only conditional on the exposure rate to the disease.
In this article, we integrate different aspects of resource allocation and create a synergistic intervention strategy that gives vulnerable populations higher priority in medical resource distribution. 
This prevention-centered strategy seeks a balance between geographical coverage and social group fairness. 
The proposed principle can be applied to other scarce resources and social benefits allocation.
\end{abstract}

\section{Introduction}\label{sec:intro}

The COVID-19 pandemic has been widely spreading around the globe and caused hundreds of thousands of deaths, crashed the healthcare systems of many countries, and stalled almost all social and economical activities, which leads to an astronomical amount of financial loss. 
To battle the spreading COVID-19 pandemic, several leading countries and organizations have devoted a significant amount of resources to develop vaccines, new diagnostics, and anti-infective treatments for the novel coronavirus \cite{whocandid}. More than 60 candidate vaccines are now in development worldwide, and several have entered early clinical trials in human volunteers, according to the World Health Organization \cite{le2020covid}. Once the development of the treatments is approved for public use, 
it is of paramount importance that these resources are quickly dispatched to all the communities in a country because one weak link in the defense against the virus would leave the neighboring communities even the whole state vulnerable and exposed to the spreading of the disease. However,
the existing healthcare infrastructures in many states and cities are insufficient to provide universal supplies of vaccinations and treatments. 
Hence, targeting the high risk population and setting fair principles for prioritizing the allocation accordingly will save lives, curb the further spreading of the virus, and prevent sequential global pandemic outbreaks in the following years.

In United States, due to the disparities \cite{nesbitt2016increasing} between different population subgroups in terms of health, financial stability, and accessibility to health care services, certain low-income and minority populated communities are particularly vulnerable to COVID-19.
For example, among Chicago's minority groups, the Latinx and Black groups have been hit much harder by the COVID-19 pandemic than any other ethnic groups. 
By the latest report (September 15, 2020) from the Chicago Department of Public Health (CDPH), among the death tolls caused by COVID-19, 33.0\% are from the Latinx group and 42.8\% are from Black \& non-Latinx group \cite{CDPH}. 
These two ethnic groups also lead in the number of COVID-19 confirmed cases by a large margin compared to the rest ethnic groups. 
These ratios are disproportional to the percentages of these ethnic groups of the total population in Chicago, which mainly consists of 32.3\% white, 28.7\% Hispanic, and 30.9\% African American, according to US Census Bureau\cite{statisticatlas}. 

The development of vaccine and other effective treatment medicines against any novel infectious disease generally takes a long time due to the high uncertainties in all the stages in the process, including clinical trial. 
Even if the vaccines and treatments are successfully developed, they are only available in limited quantities initially due to manufacturing, logistic and financial constraints. 
Due to the huge gap between available medical resources and the entire population in need of them, especially during the early production stages, the allocation of such resources becomes a difficult and yet pressing issue.
Moreover, due to demographic distribution and occupational factors of the residents in different geographical regions, resource allocation based on solely maximum geographical uniformity can lead to an \emph{unfair} distribution system and propagate biases across population subgroups of the population.

We design a \emph{fairness-aware} allocation framework for vaccine and scarce treatment resources considering both \emph{Geographical Diversity} and \emph{Social Group Fairness} as the guiding principles for prevention-centered strategies.
The social group fairness is based on the general fairness notion of \emph{Equality of Opportunity} \cite{hardt2016equality,elzayn2019fair}. 
An allocation strategy is fair if the average amount of resources an individual receives only depends on the individual's exposure rate to the disease and is independent of the individual's demographic or social-economic background. 
We also consider an allocation strategy to be diverse if the geographical location does not affect the averaged resources an individual can receive. 
Based on such notions of fairness and diversity, we formally define and formulate them into inequality constraints. 
% The notions of diversity and fairness are functions of the exposure probabilities, which depends on demographic distribution of populations in each region and the spread of the disease. 
%Note that the demographic distribution of population is different from an individual's demographic information. 
Our proposed research provides a solution to distribute the scarce medical resources in a fair manner to all communities and protect certain minorities and low-income groups that are more vulnerable to pandemic's impact.  
% Not only such a planning solution can help stop the spreading pandemic more effectively, but also pushes for justice and fairness in healthcare decision making. 
Not only can such a solution help stop the spreading pandemic more effectively, but also push for justice and fairness in healthcare decision making. 

% Our proposed research provides a solution to distribute the scarce medical resources in a fair manner to all communities, and particularly protect certain minorities and low-income groups that turn out to be vulnerable to this pandemic.  
% Not only such a planning solution can help stop the spreading pandemic more effectively, but also pushes for justice and fairness in the healthcare decision making. 

\subsection{Related Work}\label{sec:related}

Resource allocation has been a classic and important problem in many domains such as economy (e.g,\cite{monopoly}), management science (e.g,\cite{managing}), emergency response (e.g,\cite{emergency}), etc. \new{
More recent works have addressed the need for an integrated approach with different purposes such as material allocation in a production plan \cite{Added-paper1}, supplier selection, and order allocation application\cite{Added-paper2}.}
Healthcare resource allocation is one of the most challenging allocation problems. 
A large number works have been developed to identify effective strategies for it \cite{becker1997optimal,tanner2008finding,yarmand2014optimal,matrajt2013optimal,feng2017multiobjective,ho2019branch,bertsimas2013fairness,emanuel2020fair}.

Recently, the fairness of algorithmic decision making has been the center of attention of many researchers, who are designing and developing algorithms for different purposes, such as machine learning, ranking, and social welfare \cite{barocas2016big,fairmlbook,vzliobaite2017measuring,dubey2020representation}. 
Mitigating the bias of an outcome from a decision model, which is mainly caused by the inherent bias in the data and societal norm, will ensure that the outcome is not favorable or adversarial toward any specific subgroups \cite{vzliobaite2017measuring,fairmlbook,corbett2017algorithmic,dwork2012fairness,hardt2016equality,zehlike2017fa,yao2017beyond}. 
One of the critical problems where the fairness of the outcome matters significantly is scarce resource allocation.
The notion of fairness in resource allocation has been introduced in \cite{fairness_old} and later with a more precise definition in \cite{price}. 

% Recently, fairness of algorithmic decision making has been the center of attention of many researchers designing and developing algorithms for different purposes, such as machine learning, ranking, and social welfare \cite{barocas2016big,fairmlbook,vzliobaite2017measuring,dubey2020representation}. Mitigating the bias of an outcome from a decision model, which is mainly caused by the inherent bias in the data and societal norm, will ensure that the outcome is not favorable or adversarial toward any specific subgroup of observations \cite{vzliobaite2017measuring,fairmlbook,corbett2017algorithmic,dwork2012fairness,hardt2016equality,zehlike2017fa,yao2017beyond}. One of the critical problems where fairness of the outcome matters significantly is scarce resource allocation.
% The notion of fairness in resource allocation has been introduced in \cite{fairness_old} and later with a more precise definition in \cite{price}. 
%Fairness in resource allocation has been addressed in many applications such as computer networks and bandwidth allocation \cite{lan2010axiomatic,huaizhou2013fairness,bertsimas2012efficiency}, education\cite{byrne1974planning,hnat2015distributive}, 

Fairness was firstly adopted in bandwidth allocation problem for computer network systems \cite{chiu1989analysis,demers1989analysis,kelly1998rate,jaffe1981bottleneck,ogryczak2007multicriteria,ogryczak2005telecommunications,lan2010axiomatic,huaizhou2013fairness,bertsimas2012efficiency}. 
In these settings, the amount of resources requested can be modified by different users. 
Besides, service allocation mainly covers a group of users and is not necessarily one-to-one allocation. 
These settings differ from the resource allocation for scarce treatments which is a one-to-one allocation problem.

Now fairness has been addressed in many resource and service allocation methods \cite{verweij2009moral, matrajt2013optimal,lum2016predict, elzayn2019fair, koonin2020strategies}. \new{The importance of ethical consideration in resource allocation and principal guidelines have been discussed in \cite{yip2021healthcare}.} 
Here we highlight a few and emphasize our difference from them. 
In \cite{elzayn2019fair}, the authors formalize a general notion of fairness for allocation problems and investigate its algorithmic consequences when the decision-maker does not know the distribution of different subgroups (defined by creditworthy or criminal background) in the population. The distribution estimation is accomplished using censored feedback (individuals who received the resource, not the true number).
In our work, we estimate the distribution of different social groups from the data using Bayes rule. 
Singh \cite{singh2020fairness} considers a fair allocation of multiple resources to multiple users and have proposed a general optimization model to study the allocation. 
Our proposed method differs from \cite{singh2020fairness} on two aspects. 
First, \cite{singh2020fairness} assumes multi-resources and multi-type users, whereas we mainly focus on resource allocation problem across different regions and different population subgroups. 
Second, \cite{singh2020fairness} aims to maximize the coverage, whereas we consider a \emph{fairness-diversity} trade-off by minimizing the diversity and fairness gaps, simultaneously, across different regions. 
Donahue et al. \cite{donahue2020fairness} considers the problem of maximizing resource utilization when the demands for the resource are distributed across multiple groups and drawn from probability distributions. 
They require equal probabilities of receiving the resource across different groups to satisfy fairness and provide upper bounds on the price of fairness over different probability distributions. 
In our proposed model, we utilize a similar fairness requirement while requiring diversity (population) consideration across different regions. 
Thus, our work is an intersection between fairness and diversity. 
Furthermore, we deal with the one-to-one allocation instead of the coverage problem.  
\new{One recent works related to COVID-19 resource allocation designed a vulnerability indicator for racial subgroups that can be used as guidelines for medical resource allocation \cite{ong2021covid}. The proposed model cannot identify other vulnerable subgroups that are not geographically clustered, thus not able to form spatially concentrated communities. In our paper, we directly identify subgroups' vulnerability using exposure rates estimated from COVID-19 cases and death. We focus on the available data to estimate the necessary parameters and perform an empirical study using the proposed Algorithm \ref{alg:tune}.}

\subsection{Summary of Contributions}

% \nazanin{need to mention that we create a Resource Allocation Tool for cities to distribute scare resources (e.g., vaccines) at Zipcode level, and then justify the reasons like why holistic strategies are not applicable and interferes with fairness, also mentioning that we can not force it at individual level}

\new{In this paper, we aim to design a fairness-aware allocation strategy that considers the trade-off between geographical diversity and social fairness (demographic disparity) in allocating resources. We provide a new aspect to the classic allocation problem while seeking a synergistic intervention strategy that prioritizes disadvantaged people in the distribution of scarce resources.}
% More specifically,} we consider a fair and diverse allocation of scarce resources. 
The proposed approach is different from the existing works using the maximum utilization \cite{bolton2003consumer} or maximum coverage \cite{singh2020fairness,elzayn2019fair} objective.
The nature of the treatment allocation problem is not the same as the coverage problems, since the former is a one-to-one assignment problem, whereas in the latter the resource can cover more than one user such as police or doctor allocation\cite{elzayn2019fair,donahue2020fairness}. 
In our work, we aim to study scarce resource allocation considering the trade-off between geographical diversity and social fairness. 
More specifically, we do not seek to share the vaccines equally among groups of populations. 
Instead, we aim to emphasize the protection of the vulnerable populations or the ones at high risk, 
depending on their exposure rate, in order to prevent death/spread of the disease. We focus on the available data to estimate the necessary parameters and perform an empirical study using the proposed Algorithm \ref{alg:tune}. 
However, our proposed model can be generalized to resource allocation in other scenarios such as disaster relief.

In \S~\ref{sec:technical},
% we first provide a mathematical definition of
we first model the notions of geographical diversity and social group fairness for the allocation of scarce medical treatments and vaccines.
We then formulate the allocation problem into an Integer Program (IP) problem, which incorporates the diversity and fairness as constraints, in addition to capacity constraint, i.e., the amount of available resources.  
Fairness and diversity constraints are bounded by user-defined hyperparameters, $\epsilon_d$, and $\epsilon_f$, which are the allowed diversity and fairness gaps, respectively.
To obtain a feasible solution for the original IP problem efficiently, we relax the IP problem to a Linear Programming (LP) problem. 
Moreover, the fairness and diversity constraints are much more complicated than the capacity constraint, and to deal with them efficiently, we use the penalty method, which is a common practice to solve constrained optimization. 
The two constraints are combined into a single objective function using a trade-off hyperparameter $\alpha$.
To guarantee that the converted problem is equivalent to the original feasibility problem we provide theoretical proofs and subsequently,
a binary search algorithm is used to obtain a feasible range of $\alpha$.
We evaluate the allocation scenario under different $(\epsilon_d, \epsilon_f)$ values and provide the corresponding feasible range for $\alpha$, accordingly. 
Different levels of the trade-off between diversity and fairness are presented. 
The proposed framework can be applied at different stages of the pandemic to estimate the exposure rate of population subgroups, and obtain a feasible allocation considering both population size and exposed population.
In \S~\ref{sec:exp}, we evaluate the performance of the proposed model using COVID datasets in Chicago for vaccine and scarce treatment allocations. 
The results demonstrate the impact of incorporating fairness criteria in the allocation model compared to the diverse allocation and uniform allocation.
The paper concludes in \S~\ref{sec:con}.

\section{\textcolor{black}{Problem Definition}}\label{sec:technical}
%We focus on the zero phase of the vaccine distribution at the beginning of an epidemic, when only a few doses of vaccine stockpiles are available and has not been distributed before.
%This zero phase of vaccine distribution is critical to the intervention of the outbreak.
%We also assume that the number of individuals who are naturally immuned is negligible to affect the distribution policy.
% \hadis{problem description how and why? equal geographical distribution (we need to look at the diversity definition) just look at population proportionality and miss the regions with higher risk communities. to fairly prioritize the regions,...}
% \hadis{may move the first part to intro}
\new{
Traditional resource allocation treatments have been widely studied in the optimization literature \cite{monopoly,managing,emergency} where the limited resource is being distributed to areas based on a single objective (e.g., function of total population). However, regions might have different needs or require higher priority in a medical resource allocation setting (e.g.vaccine). An allocation treatment that is solely based on geographical diversity, suggests an equalized distribution among the regions, which may not provide equity in the sense of social fairness. To assure the concept of ``equity'' as well as ``equality'', we aim to model and incorporate fairness notion on demographic subgroups as another principal in the allocation decision process. 
Identifying each region's priority or risk level, is a critical task since it can save many lives and protect the vulnerable sensitive subgroups. 
%That is,we aim to design a fair and diverse resource allocation framework through  analyzing the fairness and diversity trade-off in the optimization setting. We then propose a control hyperparameter($\alpha$) for the trade-off in order to provide the best possible allocation solution. 
% That is, we seek to design fair and diverse scarce resource allocation strategies. We develop an optimization problem to minimize \emph{Social Group Fairness} and \emph{ Geographical Diversity} 
We aim to design a fair and diverse resource allocation framework through modeling the \emph{Social Group Fairness} and \emph{Geographical Diversity} and modeling their trade-off in the optimization setting.
The optimization model minimizes the fairness and diversity gaps across different subgroups of different regions while satisfying the capacity constraint to protect a vulnerable population. A trade-off analysis is performed to demonstrate the price of fairness incorporation as the fairness gap (the maximum allocation gap between demographic subgroups) with and without fairness consideration in the allocation.
We propose a tuning approach to identify an optimal range for the trade-off hyperparameter($\alpha$) in order to provide the best possible allocation solution under different scenarios.} 
% We aim to design a fair and diverse resource allocation framework through modeling the \emph{Social Group Fairness} and \emph{ Geographical Diversity} and modeling their trade-off in the optimization setting. 

\subsection{\textcolor{black}{Notations}}

% \nazanin{We need to modify some unclear notations here.}
% \nazanin{We need to add serial numbers for better organization of equations}

We consider a centralized decision maker for allocating available $b$ units of vaccines to a set of centers denoted by $M$.
They include clinics, hospitals, pharmacies, etc.
For convenience, we assume the entire area covered by the $M$ centers to be a city, but it can be a county, or other administrative district. 
Let $x_j$ be the decision variable denoting the amount of vaccines to be allocated to center $j\in M$.
Let $z_{j}$ be the region, such as the list of zip-code areas, that are assigned to be covered by center $j$.
For simplification, we assume there is no overlap between the list of zip-code areas covered by different centers (even if they are close in distance), i.e. $z_l \cap z_k=\emptyset, \forall l, k\in M$.
We also assume the policy that residents can only receive the resources from the center that covers the region where they reside. 
Such policy is not uncommon in practice, especially in distribution of scarce resources. 

\new{Next we introduce some key concepts and their notation. 
If we consider the geographical location of any individual reside to be a random variable, denoted by $Z$, then $\{j\in M, z_j\}$ are the possible values for $Z$.
Therefore, $P(Z=z_j)$ is the proportion of the population who reside in $z_j$.
Let $U_1,\dotsc,U_p$ be discrete-valued sensitive variables corresponding to demographic and socioeconomic attributes, and $S_i$ for $i=1,\ldots,p$ be the set of possible levels for each of the sensitive variable.
For instance, if $[U_1,U_2,U_3]$ represent three attributes, income, race, gender, respectively, then $S_1=\{\text{low, medium, high}\}$, $S_2=\{\text{black, latinx, white, others}\}$, and $S_3=\{\text{female, male}\}$.
The combinations of all the levels of $U_1,\ldots, U_p$ is a set denoted by $\mathcal{G}$ and indexed by a set $I$, i.e., $\mathcal{G}=\{g_i, \text{ for }i \in I\}$. 
In other words, $\mathcal{G}=S_1\times S_2\times\cdots\times S_p$.
For any $g\in \mathcal{G}$, it corresponds to a possible combination of levels of $U_1,\ldots, U_p$, such as $<low-income, black, female>$, and there can be a group of population whose values of the sensitive variables $(U_1,\ldots, U_p)$ are equal to $g$.
For short, we call it social group $g$. 
% Let the set $\mathcal{G}$ be indexed by a set $I$, i.e., $\mathcal{G}=\{g_i, \text{ for }i \in I\}$.
Let $s_{i,j}$ be the population size of the social group $g_i$ who reside in region $z_j$, where $i\in I, j \in M$.
Therefore, $s_{i,j}/\sum_{r\in I,k\in M} s_{l,j}=P([U_1,\ldots, U_p, Z]=[g_i,z_j])$.
Let $E$ be a binary random variable with $E=1$ representing the individual is exposed to the infectious disease and $E=0$ otherwise.
So $P(E=1|g_i)$ is the exposure rate of the social group $g_i$ and $P(E=1|g_i, z_j)$ the exposure rate of the social group $g_i$ living the area of $z_j$. 
It is intuitive to assume these exposure rates depend on the social groups and regions. 
We will discuss more on some reasonable assumptions on the exposure rates and how to estimate them later.}

A key concept in resource allocation is the amount of the resource per capita. 
It is a ratio between the quantity of available resource and the size of the population who are going to receive the resource. 
Denote $V$ as the amount of resource one individual receives. 
One important assumption we make here is that $V$ follows a discrete uniform distribution, and there are three parameters involved, the amount of resource $X$, $\mathcal{Y}$ the population who are to receive the $X$ amount of resource, and size of the population $Y=card(\mathcal{Y})$, the cardinality of $\mathcal{Y}$. 
Therefore, the mean value of $V$ is $\E(V|X, \mathcal{Y}, Y)=\frac{X}{Y}$ is the resources per capita, and it varies with respect to the three parameters. 

The focus of this article is on the following problem.
Given a limited amount of resource $b$, such as vaccines, how should the decision maker allocate the amount of resource $x_j$ to each center $j$ satisfying \emph{geographical diversity} (quantified by $\mathcal{D}(x)$), and \emph{social fairness} (quantified by $\mathcal{F}(x)$).

% The focus of this article is on the following problem.
% Given a limited number of vaccine resource $b$, how should the decision maker allocate number of vaccines $x_j$ to each center $j$ satisfying \emph{geographical diversity} (quantified by $\mathcal{D}(x)$), and \emph{social fairness} (quantified by $\mathcal{F}(x)$), so as to inequalities considering the geographical diversity and group fairness.
% Next, we define geographical diversity and social group fairness, and the latter is based on the \emph{equality of opportunity} notion of fairness \cite{dwork2012fairness}.
% To quantify the geographical diversity and social fairness, we introduce $\mathcal{D}(x)$ and $\mathcal{F}(x)$ in Equations~\ref{eq:diversity} and \ref{eq:fairness}.

% To simplify the notation, we define $s_{i,.}=\sum\limits_{j\in M} s_{i,j}$ as the population size in group $g_i$, $s_{.,j}=\sum\limits_{i\in I} s_{i,j}$ as the population size in $z_j$, and $S=\sum_{i,j} s_{i,j}$ is the total population.
% We also assume all the available resources are distributed to all the centers, i.e., $\sum\limits_{j\in M} x_j=b$.
\subsection{\textcolor{black}{Diversity and Fairness Modeling}} 
In this part, we define geographical diversity and social group fairness, which the latter is based on the \emph{equality of opportunity} notion of fairness \cite{dwork2012fairness}.
To quantify the geographical diversity and social fairness, we introduce $\mathcal{D}(x)$ and $\mathcal{F}(x)$ in Equations \eqref{eq:diversity} and \eqref{eq:fairness}.

\begin{definition}{\textbf{(Geographical Diversity)}}\\
% \nazanin{need to define $s_{i,j}$}
Geographical Diversity of allocation of limited resource to a set of centers $M$ is satisfied if $\forall j \in M$, on average the resource per capita is invariant with respect to the location of the groups of the population, i.e., $\E(V|x_j, z_j, \sum_{i\in I}s_{i,j})$ does not vary with respect to the location $z_j$.
\end{definition}
The notation $\E(V|x_j, z_j, \sum_{i\in I}s_{i,j})$ is simplified and we use $z_j$ to refer to the population who reside in region $z_j$. 
Let \new{$s_{i,j}$ be the population size of the social group $g_i$} and $\E(V|\sum_{j\in M} x_j, \sum_{i\in I,j\in M}s_{i,j})$ denote the averaged resource per capita over the entire city under the consideration of the resource allocation plan, and the location is omitted since it is obvious. 
It is straightforward to formulate the geographical diversity $\mathcal{D}_j(x)$  for $\forall j \in M$ as follows.
\begin{equation}\label{eq:diversity}
\mathcal{D}_j(x) = \abs*{\E(V|x_j, z_j, \sum_{i\in I}s_{i,j})-\E(V|\sum_{j\in M}x_j, \sum_{i\in I,j\in M}s_{i,j})}=\abs*{\frac{x_j}{\sum\limits_{i\in I}s_{i,j}}-\frac{\sum\limits_{j\in M} x_j}{\sum\limits_{j\in M} \sum\limits_{i\in I} s_{i,j}}}.
\end{equation}
So if geographical diversity is strictly met, $\mathcal{D}_j(x)=0$ for all $j$.

The geometrical diversity represents the conventional or minimal requirement for resource distribution, that the resources should be evenly distributed among all geographical regions across the entire city. 
Next, we introduce the concept of social fairness. 
In this definition, we emphasize the even distribution of resource among the endangered population, i.e., who are exposed to the infectious disease, disregarding the social groups. 

\begin{definition}{\textbf{(Social Fairness)}}\\
% \nazanin{need to define $s_{i,j}$}
Social Fairness of allocation of limited resource to a set of centers $M$ is satisfied if $\forall i \in I $, the averaged resource per capita is invariant with respect to the values of $U_1,\dotsc,U_p$ of the group of population, i.e., $\E(V|E=1, g_i)$ does not vary with respect to the social group $g_i$. 
\end{definition}

The notation $\E(V|E=1, g_i)$ is simplified and we use $E=1$ and $g_i$ to denote the exposed individuals in social group $g_i$.
Directly translating this definition into formula, the fairness principle should be
\begin{equation}
     \mathcal{F}_i(x)=\abs*{\E(V|E=1, g_i) - \E(V|E=1) }, \quad \forall  i\in I. 
\end{equation}

However, the calculation of $\E(V|E=1, g_i)$ and $\E(V|E=1)$ are not so straightforward, and we derive them as follows. 

We first calculate the resource per capita for the exposed individuals who reside in $z_j$, disregarding the social group, i.e., 
\begin{equation}
    \E(V|E=1, z_j)= \frac{x_j}{\sum\limits_{l \in I} s_{l,j}\times P(E=1|z_j, g_l)}, \quad  \forall j \in M
\end{equation}

in which the denominator is the amount of exposed individuals in $z_j$. 
Then,
\begin{subequations}
\begin{align}
& \E(V|E=1, g_i)=\sum\limits_{j \in M} \E(V|E=1, z_j)P(Z=z_j|E=1, g_i)\\
=&\sum\limits_{j\in M} \E(V|E=1, z_j)\frac{P(E=1, [U_1,\ldots, U_p]=g_i, Z=z_j)}{P(E=1, [U_1,\ldots, U_p]=g_i)}\\
=&\sum\limits_{j\in M} \frac{x_j}{\sum\limits_{l \in I} s_{l,j}\times P(E=1|z_j, g_l)}\frac{s_{i,j}P(E=1|g_i,z_j)/(\sum_{i', j'} s_{i', j'})}{\sum\limits_{k\in M} s_{i,k}P(E=1|g_i, z_k)/(\sum_{i', j'} s_{i', j'})}\\
=&\sum\limits_{j\in M} \frac{x_j}{\sum\limits_{l \in I} s_{l,j}\times P(E=1|z_j, g_l)}\frac{s_{i,j}P(E=1|g_i,z_j)}{\sum\limits_{k\in M} s_{i,k}P(E=1|g_i, z_k)}
\end{align}
\end{subequations}

\begin{figure}
\centering
% \begin{subfigure}{0.32\textwidth}
% \includegraphics[width=\linewidth]{plots/White_Chicago.png}
% \caption{\small White\label{fig:white}}
% \end{subfigure}
% \begin{subfigure}{0.32\textwidth}
% \includegraphics[width=\linewidth]{plots/Hispanic_Chicago.png}
% \caption{\small Hispanic groups.\label{fig:hispanic}}
% \end{subfigure}
% \begin{subfigure}{0.32\textwidth}
% \includegraphics[width=\linewidth]{plots/Black_Chicago.png}
% \caption{\small African-American.\label{fig:black}}
% \end{subfigure}
\subfigure[White]{\label{fig:white}\includegraphics[width=40mm]{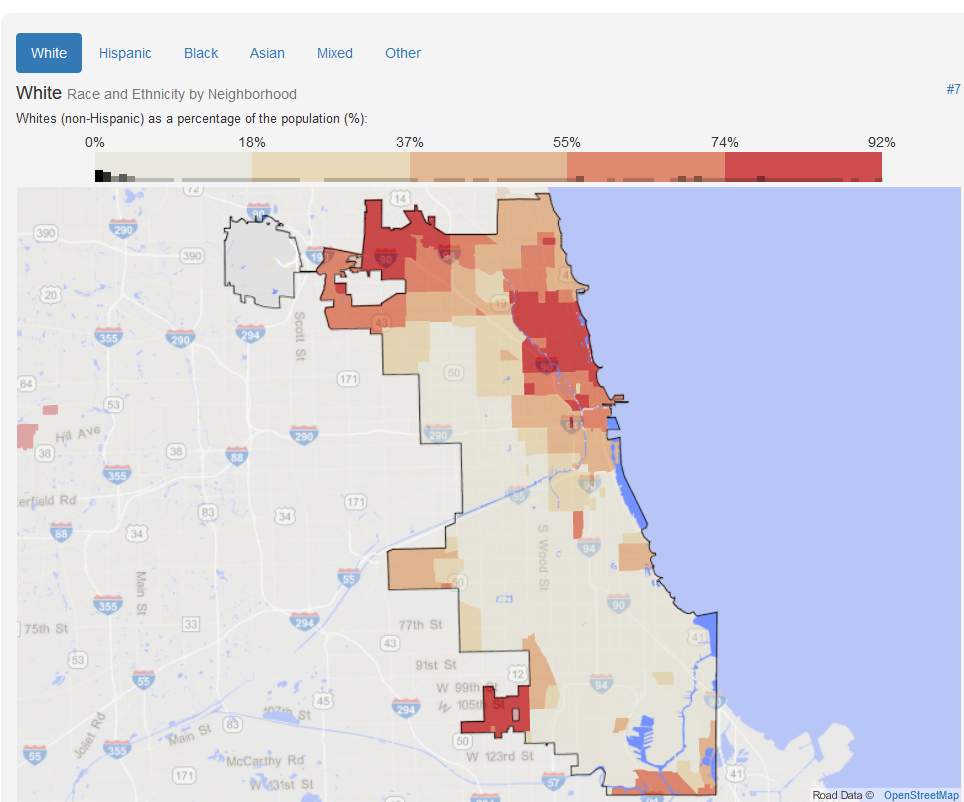}}
\subfigure[Hispanic groups.]{\label{fig:hispanic}\includegraphics[width=40mm]{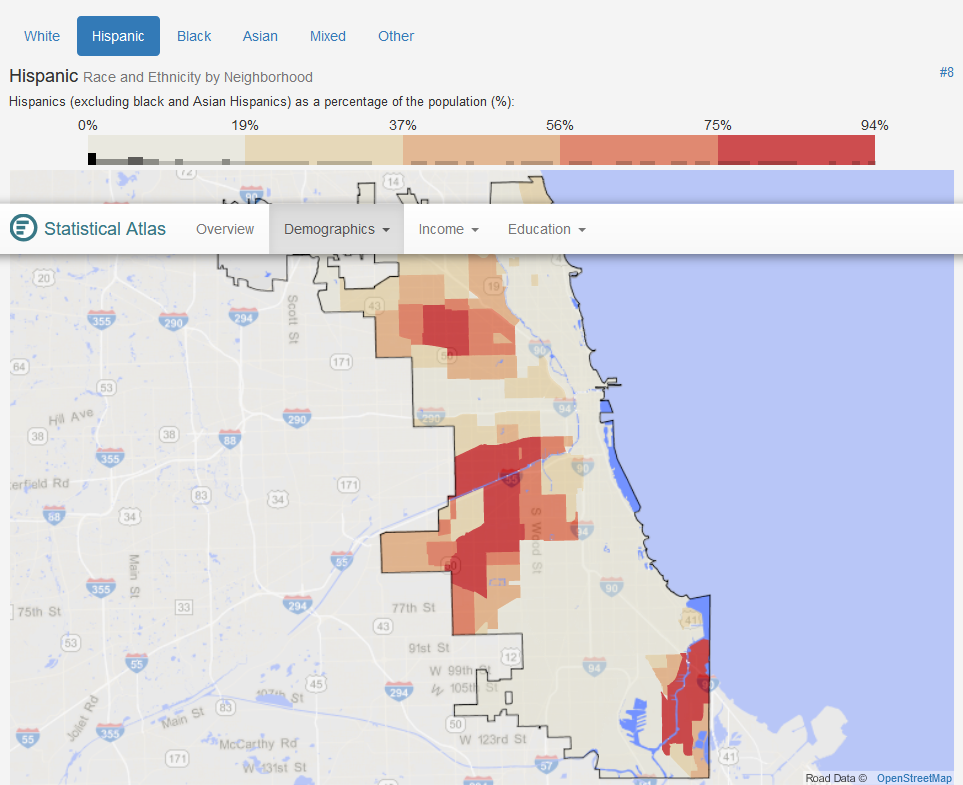}}
\subfigure[African-American]{\label{fig:black}\includegraphics[width=40mm]{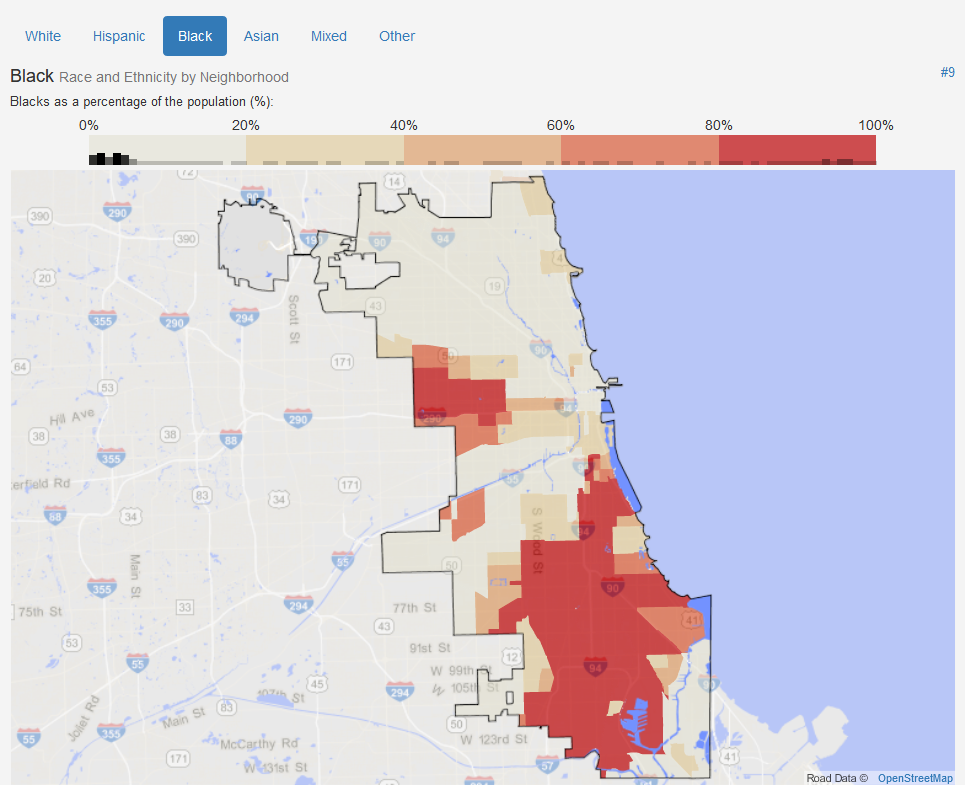}}
\caption{Map of Race and Ethnicity by Neighborhood in Chicago.\label{fig:racemap}}
\end{figure}

In this article, we assume that the chance of exposure of an individual only depends on the individual's social attributes, and is independent of the geographical location, i.e., $P(E=1|U_1,\ldots, U_p, Z)=P(E=1|U_1,\ldots, U_p)$.
Equivalently, $\forall i \in I$ and $\forall j \in M$, $P(E=1|g_i, z_j)=P(E=1|g_i)$.
This assumption is reasonable since in many U.S. cities as in Chicago, the geometrical locations of the residents are in fact highly correlated with the social and economical status of the residents.
For example, in Figure \ref{fig:racemap}, the percentage of the major ethnic groups, White, Hispanic, and African-American, are shown in three heat maps. 
It is very clear that the geometrical locations of the residents and the racial groups are heavily correlated. 
Of course, this assumption also makes the rest of the formulation much simpler. 
Under this assumption, we can simply obtain
\begin{equation}\label{eq:E=1_g_i}
\E(V|E=1, g_i)=\sum\limits_{j \in M}\frac{x_j}{\sum\limits_{l \in I} s_{l,j}\times P(E=1|g_l)}\frac{s_{i,j}}{\sum\limits_{k\in M} s_{i,k}}.
\end{equation}

Next, to obtain $\E(V|E=1)$, we need to integrate Equation \eqref{eq:E=1_g_i} with respect to $g_i$, i.e.,
\begin{subequations}
\begin{align*}
\E(V|E=1)&=\sum\limits_{i \in I} \E(V|E=1, g_i)P([U_1,\ldots, U_p]=g_i|E=1)\\
&=\sum\limits_{i \in I} \E(V|E=1, g_i)P(E=1|g_i)\frac{P([U_1,\ldots, U_p]=g_i)}{P(E=1)}\\
&=\sum\limits_{i \in I} \E(V|E=1, g_i)P(E=1|g_i)\frac{\sum\limits_{k\in M} s_{i,k}/(\sum_{i', j'} s_{i', j'})}{\sum\limits_{r\in I, k'\in M} s_{r,k'} P(E=1|z_k', g_r)/(\sum_{i', j'} s_{i', j'})}\\
&=\sum\limits_{i \in I} \left(\sum\limits_{j \in M}\frac{x_j}{\sum\limits_{l \in I} s_{l,j}\times P(E=1|g_l)}\frac{s_{i,j}}{\sum\limits_{k\in M} s_{i,k}} \right)\frac{P(E=1|g_i) \sum\limits_{k\in M} s_{i,k}}{\sum\limits_{r\in I, k'\in M} s_{r,k'} P(E=1|g_r)}\\
&=\sum\limits_{i \in I} \sum\limits_{j \in M}\frac{x_j}{\sum\limits_{l \in I} s_{l,j}\times P(E=1|g_l)}\frac{P(E=1|g_i)s_{i,j}}{\sum\limits_{r\in I, k'\in M} s_{r,k'} P(E=1|g_r)}.
\end{align*}
\end{subequations}

Based on the derivations, we can formulate the fairness for $\forall i\in I$, 
\begin{align}\label{eq:fairness}
& \mathcal{F}_i(x)=\abs*{\E(V|E=1, g_i) - \E(V|E=1) } = \\\nonumber
& \abs*{
\sum\limits_{j \in M}\frac{x_j}{\sum\limits_{l \in I} s_{l,j}\times P(E=1|g_l)}\frac{s_{i,j}}{\sum\limits_{k\in M} s_{i,k}}-
\sum\limits_{i \in I} \sum\limits_{j \in M}\frac{x_j}{\sum\limits_{l \in I} s_{l,j}\times P(E=1|g_l)}\frac{s_{i,j}P(E=1|g_i)}{\sum\limits_{r\in I, k'\in M} s_{r,k'} P(E=1|g_r)}
}. 
\end{align}
Similar to the definition of $\mathcal{D}_j(x)$, if the fairness is strictly met, $\mathcal{F}_i(x)=0$ for all $i\in I$.

% First, $P(V=1|E)$ is the probability of an individual receiving the resource, disregarding his/hers geographic and social information.
% \begin{align}
% & \mathcal{F}_i(x) = \abs*{P(V=1|E=1, g_i)-P(V=1|E=1)}= \\\nonumber\label{eq:fairness2}
% & \abs*{\sum\limits_{j \in M}\frac{ x_j }{\sum\limits_{i\in I} s_{i,j} \times P(E=1|g_i)}\times\frac{s_{i,j}}{\sum\limits_{k\in M} s_{i,k}}-\sum\limits_{i \in I} \frac{\sum\limits_{j\in M}s_{i,j}}{\sum\limits_{k\in I} \sum\limits_{j\in J}s_{k,j}}\sum\limits_{j \in M}\frac{ x_j }{\sum\limits_{i\in I} s_{i,j} \times P(E|g_i)}\times\frac{s_{i,j}}{\sum\limits_{k\in M} s_{i,k}}}, \quad \forall i \in \mathcal{I}.
% \end{align}
%\hadis{I updated the changes in tech section up to here}
%Similar to the definition of $\mathcal{D}_j(x)$, if the fairness is strictly met, $\mathcal{F}_i(x)=0$ for all $i\in I$. 

\section{\textcolor{black}{Fair and Diverse Allocation Optimization}}\label{sec:opt}

%More particularly, $\mathcal{D}_j(x)$ refers to the absolute difference between the fraction of the population in region $j$ to be
%vaccinated and the overall fraction of the population to be vaccinated. Similarly,  $\mathcal{F}_i(x)$ refers to the absolute difference between the fraction of the exposed population $i$ in region $j$ to be
%vaccinated and the overall fraction of the exposed population to be vaccinated. %The first component in Equation~\ref{eq:fairness}, corresponds to the $\frac{\text{probability of an individual in group $i$ being vaccinated}}{\text{probability of group $i$ being exposed to the disease}}$.

% We aim to minimize the diversity and unfairness gaps and ideally equal to zero to satisfy the independence definitions. To do so, we consider the maximum discrimination of allocation in geographical regions, $\mathcal{D}(x)=max_{j\in M}\mathcal{D}_j(x)$, and population subgroups, $\mathcal{F}=\max_{i\in I}\mathcal{F}(x)$, to be upper bounded by $\epsilon_{\mathcal{D}}$ and $\epsilon_{\mathcal{F}}$, respectively.
% In particular, we want to allocate our $b$ available doses of vaccines among different population subgroups in different regions. 

As explained above, 
if the diversity and fairness requirements are strictly met, all $\mathcal{D}_j(x)$ and $\mathcal{F}_i(x)$ should be 0. 
However, such constraints are too restrictive and can be difficult or impossible to satisfy for all regions and social groups. 
Define $\mathcal{D}(x)=\max_{j\in M}\mathcal{D}_j(x)$ and $\mathcal{F}(x)=\max_{i\in I}\mathcal{F}_i(x)$. 
$\mathcal{D}(x)$ and $\mathcal{F}(x)$ are auxiliary decision variables that signify the tight upper bounds on the diversity and fairness constraints, i.e., $\mathcal{D}_j(x)\leq \mathcal{D}(x)$ and $\mathcal{F}_i(x)\leq \mathcal{F}(x)$ for any $j \in M$ and $i \in I$. 
% They should be minimized since we want to achieve the diversity and fairness as much as possible. 
Ideally, we want to find the a feasible solution $x$ such that both upper bounds are equal to zero. 
In addition, $x_j$ for $j=1,\ldots, M$ should satisfy the capacity constraint, i.e., $\sum_{j\in M} x_j=b$. 
% Not considering other (financial or timeline) costs in this resource allocation problem and only 
Seeking to achieve the geographical diversity and social fairness simultaneously, 
one can formulate this problem as multi-objective (MO) minimization.

%We wish to find the vector $X\geq 0$ that satisfies these constraints considering the capacity constraint. Let us consider the feasibility problem of the following form:
% Instead of strict equality of diversity and fairness as a result of the independence of the corresponding probabilities from the geographical region and population subgroups,
%Note that, our problem formulation is agnostic to the fairness definition and many other previously studied definitions of fairness can also be used
%as $\mathcal{F}(X)$.
%\hadis{this is group fairness explain}
%\hadis{The goal is to achieve geographical diversity and fairness!}
% If we do not consider other financial or time costs in this resource allocation problem and only seek to achieve the geographical diversity and social fairness, we can formulate the following feasibility problem. 

\begin{align*}
P1: \quad \min_x \quad & (\mathcal{D}(x), \mathcal{F}(x))\\
&\mathcal{D}_j(x)\leq \mathcal{D}(x), \quad \forall j \in M\\
&\mathcal{F}_i(x)\leq \mathcal{F}(x), \quad \forall i \in I\\
&\sum_{j\in M} x_j= b\\
&\;x_j \geq 0, \quad\quad \forall j \in M\\
& x_j\in \mathbb{Z}, \quad\quad \forall j \in M.
\end{align*}

The integer constraint is because usually, resource such as vaccines are counted in integers and one individual only receives one vaccination. The integer constraints in $P1$ can be relaxed, particularly in practice when $b$ is large. 

%We relax this constraint in computation as shown later.  

%Assuming that the objective function of \hadis{this needs to be revised}$P1$ is a constant value,
%scalarization approach
% \begin{align}
% P1:\\%\label{eq:objfunc} & \max 0^T X \\
% %\nonumber &\mbox{ s.t. }\\
% \label{eq:const1} &\mathcal{D}_{j}(x)<\epsilon_d \quad\quad \forall j \in M\\
% \label{eq:const2} &\mathcal{F}_{i}(x)<\epsilon_f  \quad\quad \forall i \in I\\
% \label{eq:const3} &\sum_{j\in M} x_j= b\\
% \label{eq:boxconst} &\;x_j \geq 0, \quad\quad \forall j \in M\\
% \label{eq:int} &\;x_j\in \mathbb{Z}, \quad\quad \forall j \in M
% \end{align}

% \hadis{Since the diversity and fairness measures are defined as absolute deviations, Equation~\ref{eq:fairness} and Equation~\ref{eq:diversity}, and $\mathcal{D}(x)$ and $\mathcal{F}(x)$ are the maximum of these deviations, the feasible region in $P1$ is non-convex. Hence, inspired by Lagrangian logic, we combine the diversity and fairness constraints with a trade-off control hyperparameter $\alpha$, and convert the problem to a single objective convex optimization problem as follows:}

%\hadis{do we need the constraints Dj < D?}
A common solution to the multi-objective optimization problem is to use the weighted sum method, which leads to a simpler minimization problem $P2$. After removing the absolute operation in the constraints, we obtain the relaxed linear programming (LP) problem.

% \begin{align*}
% P2:\quad \min_{x}\quad  & (1-\alpha) \mathcal{D}+\alpha \mathcal{F}\\ 
% \mbox{ s.t. }  &\mathcal{D}_j(x)\leq \mathcal{D}, \quad \forall j \in M\\
% &\mathcal{F}_i(x)\leq \mathcal{F}, \quad \forall i \in I\\
% &\sum_{j\in M} x_j= b\\
% &x_j \geq 0, \quad\quad \forall j \in M\\
% &x_j\in \mathbb{Z},\quad\quad \forall j \in M.
% \end{align*}

% where $\alpha \in [0,1]$ is the control parameter on the importance of fairness for allocation. If $\alpha >0.5$ the objective function focuses more on the fairness in the minimization and $\mathcal{F}(X)$ will be smaller more likely. Meaning that $\mathcal{F}$ deviates less from the constraint ~\ref{eq:const2}.
\begin{align*}
P2: \quad  \min_x & \quad (1-\alpha) \mathcal{D}(x) +\alpha \mathcal{F}(x)\\
\mbox{ s.t. } \quad &\mathcal{D}^+_j(x)\leq \mathcal{D}(x), \quad \forall j \in M\\
 & \mathcal{D}^-_j(x )\geq \mathcal{D}(x), \quad \forall j \in M\\
&\mathcal{F}^+_i(x)\leq \mathcal{F}(x), \quad \forall i \in I\\
& \mathcal{F}^-_i(x)\geq \mathcal{F}(x), \quad \forall i \in I\\
&\sum_{j\in M} x_j= b\\
& x_j \geq 0, \quad \forall j \in M.\\
% & x_j\in \mathbb{Z}, \quad\quad \forall j \in M.
\end{align*}

where $\mathcal{D}^+_j(x)$ refers to the positive side of the absolute value function, $\mathcal{D}_j(x)$, and $\mathcal{D}^-_j(x)$ refers to the negative side of the absolute function. Similarly, for $\mathcal{F}_j(x)$, we have $\mathcal{F}^+_j(x)$ and $\mathcal{F}^-_j(x)$.
Here $\alpha \in [0,1]$ is a hyperparameter that controls the trade-off between the importance of fairness and diversity. 
If $\alpha >0.5$, the objective function focuses more on achieving the fairness, and it focuses more on the diversity if $\alpha < 0.5$. 

\begin{wrapfigure}{R}{0.45\textwidth}
    \begin{minipage}{0.45\textwidth}
      \begin{algorithm}[H]
        \caption{{\bf Rounding Heuristic}}
        \begin{algorithmic}[1]
        \label{alg:rounding}
        \small
        \STATE Find the LP relaxation solution; solve $P3$.
        \STATE Round $x_j$ to the nearest integer for $j\in M$.
        \STATE $b'=b-\sum_{j\in M}x_j$.
        \IF{$b'<0$}
            \STATE Sort $x_j$'s according to the population size of $z_j$ (from largest to smallest).
            \STATE Decrease the $x_j$'s corresponding to the top $[b']$ populated regions by 1.
        \ELSIF{$b'> 0$}
            \STATE Sort $x_j$'s according to the exposed population of $z_j$ (from largest to smallest). 
            \STATE Increase the $x_j$'s corresponding to the top $[b']$ regions by 1.
        \ENDIF
        \STATE Return $x_j, \forall j\in M$.
        \end{algorithmic}
      \end{algorithm}
    \end{minipage}
  \end{wrapfigure}

If we relax the integer constraints of $P1$, based on the multi-objective optimization theory, it is easy to conclude that for any $x^*$ in the \emph{Pareto front} of $P1$, there exists an $\alpha^*\in [0,1]$ such that $x^*$ is an optimal solution of $P2$. 
This is because the remaining constraints of $P1$, including $\mathcal{D}_j(x) \leq \mathcal{D}(x) $ for all $j \in M$, $\mathcal{F}_i(x)\leq \mathcal{F}(x)$ for all $i\in I$, $\sum_j x_j =b$, and $x_j \geq 0$ for all $j \in M$, form a convex polyhedron. 
Particularly, the constraints $\mathcal{D}_j(x) \leq \mathcal{D}(x)$ and $\mathcal{F}_i(x)\leq \mathcal{F}(x)$ are all linear in $x$.
Thus, we have the above conclusion. 

In practice, the simplex method can be used to find the optimal solution of $P2$ efficiently with common optimization libraries like IBM CPLEX and SciPy.

Once the optimal solution of the LP-relaxation problem, $P2$, is obtained for a given $\alpha$ value, we need to round the solution into integers. 
But the rounded solution is not necessarily feasible for the minimization problem $P2$. 
To ensure the feasibility of the rounded solution we employ a heuristic rounding approach in Algorithm~\ref{alg:rounding}, which is similar to the one in \cite{tayfur2009model}.
The algorithm starts with rounding the LP relaxation solution to the nearest integer values. 
Next, if the capacity constraint $\sum_{j \in M} x_j=b$ is not satisfied, the algorithm reduces $x_j$ for the top populated region based on the total exceeded amount. 
If the resource constraint is under-satisfied, the algorithm increases the top exposed populated areas based on the total remaining resources. 
Note that in the Diverse-only ($\alpha=0$) and Fairness-only ($\alpha=1$) scenarios, the algorithm only considers the population and exposed population to sort $x_j$'s, respectively.

\subsection{Feasibility}\label{sec:feasibility}

As we discussed in \S~\ref{sec:opt}, the ideal situation would be minimizing the upper bounds $\mathbb{D}(x)$ and $\mathcal{F}(x)$ to the full extent. An important observation to make here is that the ideal case with zero upper bounds {\em both} on the fairness and the diversity constraints is almost always impossible, due to the underlying biases in data and historical discrimination in society.
Therefore, in practice, a small positive upper bound threshold is considered to be satisfied on fairness and diversity.
For example, the fairness requirement can be thought of as the US Equal
Employment Opportunity Commission’s ''four-fifths rule,''
which requires that ''the selection rate for any race, sex, or
ethnic group [must be at least] four-fifths (4/5) (or eighty
percent) of the rate for the group with the highest rate''\footnote{Uniform Guidelines on Employment Selection Procedures, 29 C.F.R. §1607.4(D) (2015).}. We consider a similar requirement for diversity as well.

\textit{The solution obtained under the MO model, $P2$, is unable to guarantee to satisfy these requirements}.
% However, the solution obtained under MO model might be either unfair or not diverse. %More particularly, minimizing the diversity and fairness gaps might prioritize one over the other due to the larger inherent gap (the one with higher gap will always dominate the other). 
For example, 
a solution with zero unfairness (but not satisfactory on the level of diversity) can be on the Pareto front solution of MO -- hence might be the optimal output -- even though it is not a valid solution.
% the Pareto front solution of MO contains diversity-only solution as well as fair-only as neither of them can be dominated by the other. 
% Since we are interested in a fair AND diverse allocation, i
% In order to achieve the {required} level of fairness and diversity at the same time, it is expected to set an upper bound for the gaps. For example, the fairness requirement can be thought of as the US Equal
% Employment Opportunity Commission’s ''four-fifths rule,''
% which requires that ''the selection rate for any race, sex, or
% ethnic group [must be at least] four-fifths (4/5) (or eighty
% percent) of the rate for the group with the highest rate''\footnote{Uniform Guidelines on Employment Selection Procedures, 29 C.F.R. §1607.4(D) (2015).}. We consider similar requirement for the diversity as well.

%Mapping $P2$ with a feasibility problem will allow us to meet the diversity and discrimination requirements for the resource allocation. %This will balance the minimization trade-off between diversity and fairness gaps.
% Suppose that to maintain diversity we need to satisfy $\mathcal{D}(x) \leq \epsilon_{\mathcal{D}}$ satisfies. 
%These control parameters are defined as hard constraints, hence, 
Subsequently, we introduce two control parameters $\epsilon_d$ and $\epsilon_f$, which are the acceptable thresholds for the diversity and fairness requirements, i.e., $\mathcal{D}(x)\leq \epsilon_d$ and $\mathcal{F}(x)\leq \epsilon_f$. 
%Both $\epsilon_d$ and $\epsilon_f$ should be relatively small, so that we can nearly achieve diversity and fairness. 
%To maintain fairness we need to satisfy $\mathcal{F}(x)\leq \epsilon_{\mathcal{F}}$.
% Our formal feasibility problem formulation is as following:
% \begin{align}
% P3:\\%\label{eq:objfunc} & \max 0^T X \\
% %\nonumber &\mbox{ s.t. }\\
% \label{eq:const1} &\mathcal{D}(x)<\epsilon_d\\
% \label{eq:const2} &\mathcal{F}(x)<\epsilon_f\\
% \label{eq:const3} &\sum_{j\in M} x_j= b\\
% \label{eq:boxconst} &\;x_j \geq 0, \quad\quad \forall j \in M\\
% \label{eq:int} &\;x_j\in \mathbb{Z}, \quad\quad \forall j \in M
% \end{align}
The violation threshold parameters $\epsilon_{\mathcal{D}}, \epsilon_{\mathcal{F}} \in [0,1]$, are user-defined values that determine how
diverse and fair the allocation should be. 
The value $\epsilon_{\mathcal{F}} = 0$ corresponds to a fully
fair allocation, whereas $\epsilon_{\mathcal{F}} = 1$ corresponds to a completely
fairness-ignorant allocation that solely considers the diversity.
These violations must be relatively small by which resources allocated to a particular region and particular group nearly achieve the required level of diversity and fairness, respectively.
% where Equation~\ref{eq:const3} is the capacity constraint; it assures that the total allocation does not exceed the capacity, i.e. what is available. Each constraints measures the deviation of each objective from a targeted goal $\epsilon$.
% $\mathcal{D}$ refers as to the geographical diversity set of constraints and $\mathcal{F}$ refers as to fairness set of constraints.
Note that as quantitative metric, we say that allocation is considered as fair if $\mathcal{F}(x)\leq \epsilon_f$, and allocation is considered as diverse if $\mathcal{D}(x) \leq \epsilon_d$. 
Depending on the constraints, the weights, and the available scarce resources, no feasible solution may exist for the optimization problem. In this case, the
decision-maker will have no choice but to relax the constraints. In \S \ref{sec:exp} we will discuss the allocation under different scenario. %\hadis{did you do talk about infeasible scenario?}

Fortunately, as we shall prove in Theorem~\ref{th:1}, if the problem has a feasible solution given the fairness and diversity constraints, there must exist an $\alpha$ value under which the optimum solution of $P2$ is feasible and vice versa.
% Theorem~\ref{th:1} shows that given the fairness and diversity constraints in, if there exist a feasible solution to $P3$ then there must exist some $\alpha$ such that the feasible solution of $P1$ is an optimum solution for $P2$.
\begin{proposition}\label{th:1}
% Let $X'\neq \emptyset$ be the set of feasible solution of $P3$. Let $A(\alpha)= (1-\alpha)\mathcal{D}+\alpha \mathcal{F})$. Given $x^* \in X'$, $\exists \alpha^* s.t. \{\argmin_x A(\alpha^*)|\sum_{j\in M} x_j\leq b, x_j \geq 0\} =x^*$.
Let be $X'=\{x|\mathcal{D}(x)<\epsilon_d, \mathcal{F}(x)<\epsilon_f, \sum_{j\in M} x_j=b, x_j \geq 0\}$:
\begin{enumerate}

    \item If $X'\neq \emptyset$: given $x^* \in X'$, $\exists \alpha^* $ such that $x^*$ is an optimal solution of $P2$.
    \item If $\not\exists ~\alpha^* $ such that $x^*$ is an optimal solution of $P2$ that satisfies both conditions $\mathcal{D}(x)<\epsilon_d$ and $ \mathcal{F}(x)<\epsilon_f$, then $X' = \emptyset$.
\end{enumerate}
% and $X'\neq \emptyset$. Given $x^* \in X'$, $\exists \alpha^* $ such that $x^*$ is an optimal solution of $P2$.
\end{proposition}
We now need to find the $\alpha$ value for which an optimum of $P2$ is a feasible satisfying diversity and fairness requirements. To do so, a naive approach would be a brute-force search in $[0,1]$ with a small step sizes added to $\alpha$. However, this is not a practical approach since it needs to solve the optimization problem in each iteration and as the step size decreases the exhaustive search increases and the number of times that it needs to solve $P2$ increases, accordingly. Therefore, we first obtain the monotone properties of optimal upper bound $\mathcal{D}(x)$ and $\mathcal{F}(x)$ with respect to $\alpha$. This result is used to find proper $\alpha$ value later in \S~\ref{sec:tuning}. 
% \begin{theorem}
% Let $x_1^*$ and $x_2^*$ be optimum solutions for $P2$ given $\alpha_1$ and $\alpha_2$, respectively, where $\alpha_2 < \alpha_1$. The monotonocity of $\mathcal{D}(x)$ and $\mathcal{F}(x)$ holds if $ \mathcal{F}(x_2^*)\leq \mathcal{F}(x_1^*)$ and $\mathcal{D}(x_2^*) \geq \mathcal{F}(x_1^*).$
% \end{theorem}
\begin{proposition}
Let $x_1^*$ and $x_2^*$ be optimum solutions of $P2$ given $\alpha_1$ and $\alpha_2$, respectively. 
It can be shown that if $\alpha_1\leq \alpha_2$, then $\mathcal{F}(x_2^*)\leq \mathcal{F}(x_1^*)$ and $\mathcal{D}(x_1^*)\leq \mathcal{D}(x_2^*)$. 
\end{proposition}
% \begin{proof}
% Let $\beta=\frac{\alpha}{(1-\alpha)}$. Given $x_1^*$ and $x_2^*$ corresponding to $\beta_1$ and $\beta_2$, where $\beta_2<\beta_1$:
% \begin{equation*}
% \mathcal{D}(x_1^*)+\beta_1\mathcal{F}(x_1^*)\leq\mathcal{D}(x_2^*)+\beta_1\mathcal{F}(x_2^*)
% \end{equation*}
% \begin{equation*}
% \mathcal{D}(x_2^*)+\beta_2\mathcal{F}(x_2^*)\leq\mathcal{D}(x_1^*)+\beta_2\mathcal{F}(x_1^*)
% \end{equation*}
% Adding the above Equations, we will have:
% \begin{equation*}
% \beta_1\mathcal{F}(x_1^*)+\beta_2\mathcal{F}(x_2^*)\leq\beta_1\mathcal{F}(x_2^*)+\beta_2\mathcal{F}(x_1^*)
% \end{equation*}
% \[
% \implies (\beta_2-\beta_1)(\mathcal{F}(x_2^*)-\mathcal{F}(x_1^*)) \leq 0
% \]
% which implies $\mathcal{F}(x_2^*)\leq\mathcal{F}(x_1^*)$. The monotonicity proof for $\mathcal{D}(x)$ is the same with $\beta=\frac{(1-\alpha)}{\alpha}$. If $\alpha_2 < \alpha_1$ then $\mathcal{D}(x_1^*)>\mathcal{D}(x_2^*)$.
% \end{proof}

\subsection{Choosing Trade-off Parameter, $\alpha$}\label{sec:tuning}

The parameter $\alpha$ in $P2$ controls the trade-off between the diversity and fairness of the allocation decision. 
Using the monotone property of $\mathcal{D}(x)$ and $\mathcal{F}(x)$ to $\alpha$ in Proposition 1, we now propose an approach to find a proper $\alpha$ value that satisfies both fairness and diversity constraints for given $(\epsilon_d,\epsilon_f)$.
%Using the monotonicity property of $\mathcal{D}(x)$ and $\mathcal{F}(x)$ in related to $\alpha$, we now propose an efficient search algorithm based on \emph{binary search} to find a valid solution for $P1$.
% \begin{wrapfigure}{L}{0.5\textwidth}
%     \includegraphics[width=\linewidth]{plots/monotone.png}
%     \vspace{-5mm}
%     \caption{\small Feasible region with respect to $\alpha$.}
%     \label{fig:monotone}
% \end{wrapfigure}

\begin{figure}[h]
    \begin{minipage}{0.5\textwidth}
        \includegraphics[width=\linewidth]{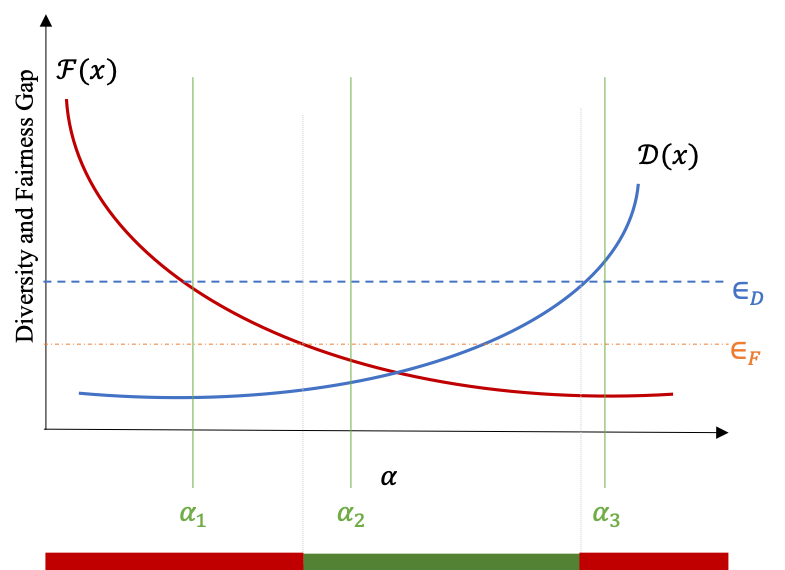}
    %\vspace{-5mm}
    \caption{\small Feasible region with respect to $\alpha$.}
    \label{fig:monotone}
    \end{minipage}
    \begin{minipage}{0.5\textwidth}
      \begin{algorithm}[H]
        \caption{{\bf Tuning $\alpha$}}
        \begin{algorithmic}[1]
        \label{alg:tune}
        \small
        \STATE $\alpha_l=0, \alpha_u=1$
            \WHILE{$|\alpha_l-\alpha_u| \geq \tau$}
                \STATE $\alpha_m=\frac{\alpha_u+\alpha_l}{2}$
                \STATE $x^*=\argmin P3(\alpha_m)$
                \IF{$\mathcal{D}(x^*) \leq \epsilon_d  ~ \& ~ \mathcal{F}(x^*)>\epsilon_f$}
                    \STATE $\alpha_l = \frac{\alpha_l+\alpha_m}{2}$
                \ELSIF{$\mathcal{D}(x^*)>\epsilon_d ~ \& ~ \mathcal{F}(x^*)\leq \epsilon_f$}
                    \STATE $\alpha_u=\frac{\alpha_u+\alpha_m}{2}$
                \ELSIF{$\mathcal{D}(x^*)\leq \epsilon_d ~\& ~\mathcal{F}(x^*)\leq \epsilon_f$}
                \STATE {\bf return} $x^*$
                \ELSE
                    \STATE {\bf return} $x^*=\emptyset$
                \ENDIF
            \ENDWHILE
            \STATE {\bf return} $x^*$
            \end{algorithmic}
        \end{algorithm}
    \end{minipage}
\end{figure}

Using an example we show the high-level idea of the tuning algorithm under different scenarios. 
In Figure~\ref{fig:monotone}, the monotonically decreasing red curve corresponds to the fairness constraint and the monotonically increasing blue curve corresponds to the diversity constraint. 
The dashed horizontal lines corresponds to the thresholds allowed for diversity $\mathcal{D}(x)$ and fairness $\mathcal{F}(x)$ constraints, respectively. 
For $\alpha < \alpha_1$ we can see that $\mathcal{F}(x)$ exceeds the allowed threshold of $\epsilon_f$. 
Based on Proposition 1, we should remove the range $[0, \alpha_1)$ from consideration. 
For $\alpha> \alpha_3$, $\mathcal{D}(x)$ exceeds the allowed threshold of $\epsilon_d$, hence, we can remove the range $(\alpha_3, 1]$ from consideration. 
For $\alpha=\alpha_2$, both $\mathcal{D}(x)\leq \epsilon_d$ and $\mathcal{F}(x)\leq \epsilon_f$ are satisfied. 
The color bar below the plot in Figure~\ref{fig:monotone} shows different ranges of $\alpha$.
The green segment corresponds to the range of $\alpha$ values for which the optimum solution of $P2$ satisfies $\mathcal{D}(x)\leq \epsilon_d$ and $\mathcal{F}(x)\leq \epsilon_f$. 
Consequently, we can prune the infeasible intervals of $\alpha$ by evaluating the violation of diversity and fairness thresholds.

% Using an example we show the high-level idea of the tuning algorithm under different scenarios. In Figure~\ref{fig:monotone}, the monotonically decreasing red curve corresponds to the fairness constraint and the monotonically increasing blue curve corresponds to the diversity constraint. The dashed horizontal lines corresponds to the thresholds allowed for diversity, $\epsilon_{\mathcal{D}}$ and fairness, $\mathcal{F}(x)$ constraints, respectively, in $P1$, Equation~\ref{eq:const1} and \ref{eq:const2}. For $\alpha=\alpha_1$ we can see that $\mathcal{F}(x)$ exceeds the allowed threshold of $\epsilon_F$. Using the monotonicity property we can prune the $\alpha$ values smaller than $\alpha_1$ because for any values in that range $\mathcal{F}(x)$ is infeasible. For $\alpha=\alpha_3$ the diversity constraint is not satisfied, hence, we can prune the range of $\alpha$ for the values larger than $\alpha_3$, for the same reason. For $\alpha=\alpha_2$, both $\mathcal{D}(x)$ and $\mathcal{F}(x)$ are satisfied. Looking at the color bar below the plot, the green part shows the range of $\alpha$ where the optimum solution of $P3$ is feasible for $P1$ and the red areas show the infeasible region of $\alpha$ values. Consequently, we can prune the $\alpha$ interval evaluating the violation of diversity and fairness constraint considering the optimum solution obtained from $P3$ according to a certain $\alpha$ value.
%\begin{wrapfigure}{L}{0.5\textwidth}
% \begin{figure}[b]
%     \begin{minipage}{0.5\textwidth}
%         \includegraphics[width=\linewidth]{plots/monotone2.png}
%     \vspace{-5mm}
%     \caption{\small Feasible region with respect to $\alpha$.}
%     \label{fig:monotone}
%     \end{minipage}
%     \begin{minipage}{0.5\textwidth}
%       \begin{algorithm}[H]
%         \caption{{\bf Tuning $\alpha$}}
%         \begin{algorithmic}[1]
%         \label{alg:tune}
%         \small
%         \STATE $\alpha_l=0, \alpha_u=1$
%             \WHILE{$|\alpha_l-\alpha_u|\geq\tau$}
%                 \STATE $\alpha_m=\frac{\alpha_u+\alpha_l}{2}$
%                 \STATE $x^*=\argmin P3(\alpha_m)$
%                 \IF{$\mathcal{D}(x^*)<\epsilon_{\mathcal{D}} \& \mathcal{F}(x^*)>\epsilon_{\mathcal{F}}$}
%                     \STATE $\alpha_l=\frac{\alpha_l+\alpha_m}{2}$
%                 \ELSIF{$\mathcal{D}(x^*)>\epsilon_{\mathcal{D}} \& \mathcal{F}(x^*)<\epsilon_{\mathcal{F}}$}
%                     \STATE $\alpha_u=\frac{\alpha_u+\alpha_m}{2}$
%                 \ELSIF{$\mathcal{D}(x^*)<\epsilon_{\mathcal{D}} \& \mathcal{F}(x^*)<\epsilon_{\mathcal{F}}$}
%                 \STATE {\bf return} $x^*$
%                 \ELSE
%                     \STATE {\bf return} $x^*=\emptyset$
%                 \ENDIF
%             \ENDWHILE
%             \STATE {\bf return} $x^*$
%             \end{algorithmic}
%         \end{algorithm}
%     \end{minipage}
% \vspace{-5mm}
% \end{figure}
%\end{wrapfigure}

Algorithm \ref{alg:tune} represents the proposed tuning algorithm. 
The algorithm starts from the middle of the $\alpha$ interval, $\alpha_m$, and solves $P2$ with that. 
Then, it splits the $\alpha$ region into half to further prune the region until a narrow feasible range of $\alpha$ is obtained. 
In each iteration, if the diversity threshold is not satisfied by the solution obtained from $P2$ the algorithm prunes the interval larger than $\alpha_m$. 
Similarly, if the fairness threshold is not satisfied, the algorithm prunes the interval smaller than $\alpha_m$. 
If neither of the thresholds is satisfied with the solution obtained from $P2$ using $\alpha_m$, there does not exist a feasible solution for $(\epsilon_d, \epsilon_f)$. 
Finally, when both constraints are satisfied the algorithm returns a valid solution for $P2$, which can be rounded to obtain an integer solution for $P1$. 

\subsection{Price of Fairness}\label{sec:pof}
% We now define the Price of Fairness (PoF) as the fairness gap that is obtained by the allocation solution without any fairness constraint to the allocation solution considering the fairness constraint.
We now define the Price of Fairness (PoF) as the difference of the fairness gap of the optimal solution of $P2$ that is obtained with and without any fairness constraints.
As we described in \S~\ref{sec:technical}, Equation~\eqref{eq:fairness} is the fairness constraint that is defined for each social group. 
Solving $P2$ with and without the fairness constraints, we can obtain different allocation solutions and consequently different fairness and diversity gaps. 
This will allow us to compare fair-diverse allocation performance with the Diverse-only allocation to analyze the impact of fairness constraints, namely PoF.

The PoF metric assists the decision-maker on the fairness scheme to be considered, mainly, $\epsilon_d$ and $\epsilon_f$. 
Note that there is a trade-off between diversity and fairness gap in $P2$, when the allocation scheme is more focused on diversity, smaller $\epsilon_f$, the control parameter $\alpha$ become close to $0$. 
Subsequently, $P2$ minimizes the diversity upper bound subject to only the capacity constraint. 
Comparing the solution of such a problem with the one that has more focus on fairness, $\alpha$ close to $1$, we expect to see a larger fairness gap and smaller diversity gap. 
In \S \ref{sec:exp}, we evaluate PoF of the allocation solution obtained under different scenarios.

\section{Case Study: Resource Allocation for COVID-19 Relief in US cities }\label{sec:exp}

In this section, we apply the proposed fair and diverse allocation strategy to the planning of distribution of medical resources for \new{COVID-19 relief among some US cities (Chicago, New York, Baltimore)}.
%Chicago's COVID-19 relief. 
\new{The city of Chicago is a segregated city with relatively high Hispanic and Black population\cite{Covid-6cities}. Several research studies have addressed a positive correlation among proportions of Black and Hispanic communities and COVID-19 cases~\cite{Covid-6cities,Covid-disparities,Covid-disparities2,Covid-disparities3}.
More specifically, in \cite{Covid-6cities}, both positive and statistically significant associations between the proportions of Black and Hispanic population and per capita COVID-19 cases  have been identified using data from six segregated cities in the US. In particular, New York City has the highest rate of confirmed cases, followed by Chicago and Baltimore. That is, we present our instance problems based on these cities.Even though our case study is limited to three cities, it includes the first (New York City) and third (Chicago) largest cities in the USA, and presents our approach for cities with higher vulnerable population (e.g., Baltimore).}

\new{We rely on the US census for population and demographic distribution data. For the COVID-19 cases, death tolls, and hospitalization we separately collect the data that is provided by governmental departments of each considered cities.
The main sources of each datasets are described in Section 3.1.} 
In the subsequent section, we first perform a descriptive analysis of the data to motivate the necessity of the fair and diverse distribution of medical resources. 
Next, we evaluate the performance of our proposed \emph{Fair-Diverse} allocation approach and identify a reasonable trade-off parameter $\alpha$ based on the binary search Algorithm \ref{alg:tune}. 
Lastly, we will calculate the price of fairness (PoF) to highlight the role of fairness constraints in our optimization setting. 
All of the analytical results, optimization models, and algorithms were implemented in Python 3.7 using docplex and sklearn packages. 
The codess are available on github \footnote{ \url{https://github.com/nnezam2/Fair-Resource-Allocation/blob/master/README.md}}.

\subsection{Data Description}\label{sec:data}

\textbf{Population Dataset:} The \emph{uszipcode}\footnote{\url{https://uszipcode.readthedocs.io/}} Python package provides detailed geographic, demographic, socio-economic, real estate, and education information at the state, city, and even zip code level for different areas within the US. Based on the documentation, this package uses an up-to-date database by having a crawler running every week to collect different data points from multiple data sources. This dataset does not provide the intersectional population. We refer to this dataset as \attrib{Pop.} throughout the experiment section.

\new{\textbf{City of Chicago COVID-19 Database}
% \textbf{Daily Chicago COVID-19 Cases, Deaths, and Hospitalizations}:  
%The \emph{City of Chicago COVID-19 database}
\cite{chicagocity} provides daily data on COVID-19 positive-tested cases, death tolls, hospitalizations, and other individuals' attributes (e.g., age, race, gender) to track the pandemic in this city. 
In our study, we primarily use the \emph{COVID-19 Daily Cases} data to reveal the inequality among different population subgroups. 
We refer to this dataset as \attrib{Chicago-COVID-Cases} throughout the experiment section. The city of Chicago COVID-19 database also provides daily data on \emph{COVID-19 Cases, Tests, and Deaths by ZIP Code} dataset \ref{sec:descriptive}. We refer to this dataset as \attrib{Chicago-COVID-Zipcode} throughout the experiment section.}

% Moreover, we access these data through the SODA API (Socrata Open Data API), which allows users to access open data resources online rather than reading from the conventional static data files.

% \textbf{COVID-19 Cases, Tests, and Deaths by ZIP Code:} The City of Chicago COVID-19 database \cite{chicagocity} also provides daily data on \emph{COVID-19 Cases, Tests, and Deaths by ZIP Code} dataset \ref{sec:descriptive}. We refer to this dataset as \attrib{COVID-Zipcode} throughout the experiment section.

% To perform an initial descriptive analysis on the COVID-19distribution across different population subgroups, we merge the \attrib{Pop.} with the \attrib{COVID-Zipcode}. We will use this dataset for descriptive analysis provided in \S \ref{sec:descriptive}. We refer to this dataset as \attrib{COVID-Pop.-Zipcode} throughout the experiment section. 

\new{\textbf{New York City (NYC) COVID-19 Data Repository\footnote{\url{ https://github.com/nychealth/coronavirus-data}}} consists of different COVID-19 related datasets including daily, weekly, monthly data, data on SARS-CoV-2 variants, the cumulative COVID-19 cases, etc. For the purpose of our analysis, we use COVID-19 cases and death totals by age, race, and gender since the start of the COVID-19 outbreak in NYC, February 29, 2020. We refer to this dataset as \attrib{NYC-COVID-Zipcode} throughout the experiment section.}

\new{\textbf{City of Baltimore COVID-19 data Dashboard\footnote{\url{https://coronavirus.baltimorecity.gov/}}} provides different statistics and visualization for COVID-19 data in Baltimore. The main source of this data is the Maryland Department of health\cite{Covid-Maryland}, which is updated on a daily basis. This dataset includes COVID-19 the number of cases and death at the zip-code level, and total cases and death tolls by age, race, gender. We refer to this dataset as \attrib{BLT-COVID-Zipcode} throughout the paper.}

\subsection{Descriptive Analysis on COVID-19 risk factor among different population subgroups, and different regions (zip-codes) }\label{sec:descriptive}

Utilizing the \attrib{Chicago-COVID-Cases} dataset, we first identify the contribution of each demographic attribute to the total number of COVID-19 cases and deaths in each region of the Chicago area. This analysis reveals key insights into the actual importance and differentiation among demographic subgroups. 

\begin{figure}[ht]
  \centering
    \begin{minipage}[b]{0.55\textwidth}
          \centering
        \includegraphics[width=1\linewidth, height=7.5cm]{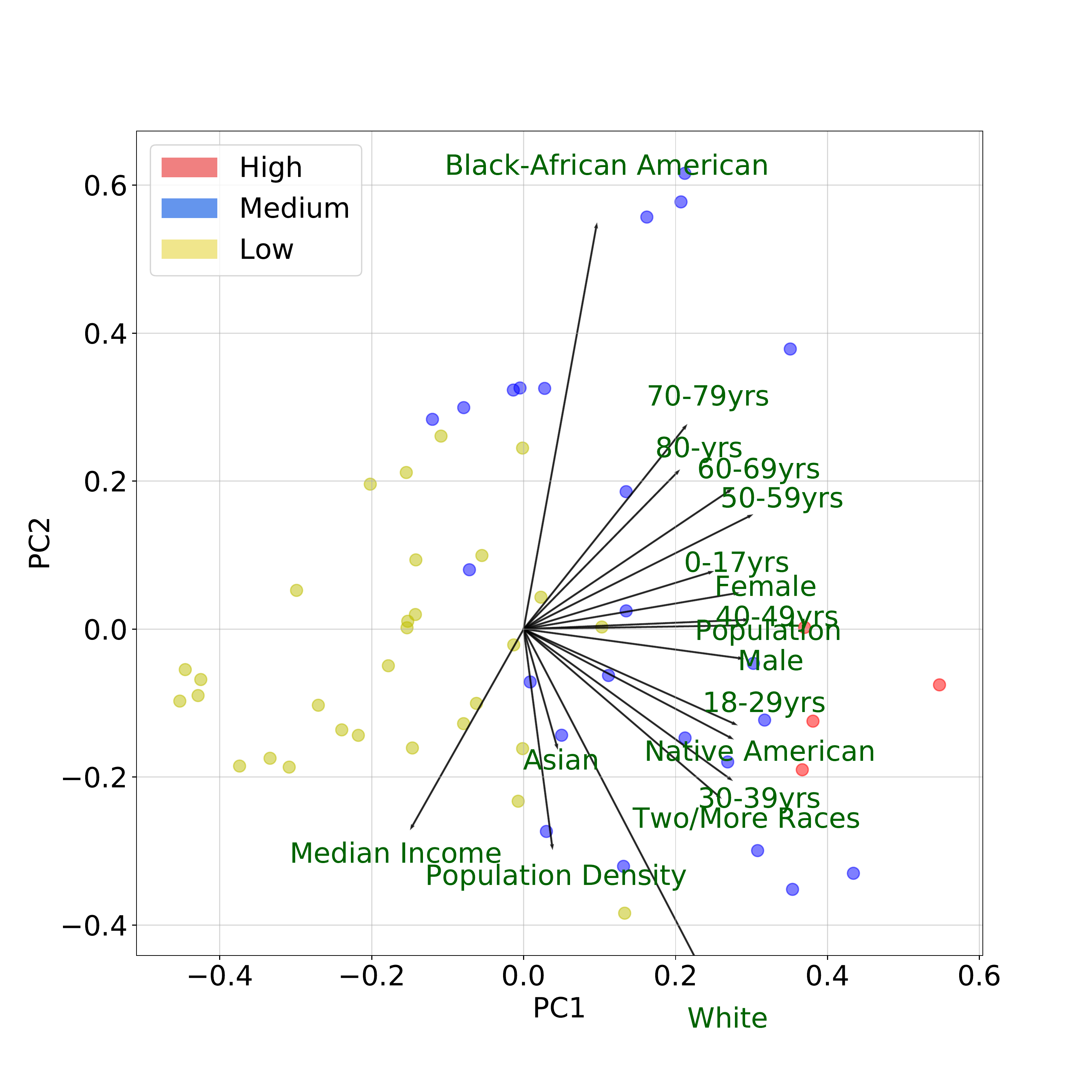}
        \caption{\new{Biplot of demographic features using PCA Analysis with Kmeans clustering -Chicago City regions}}
        \label{fig:hist}
    \end{minipage}
    \hfill
\begin{minipage}[b]{0.44\textwidth}
    \centering
  \begin{adjustbox}{width=\columnwidth,center}
  \setlength\tabcolsep{1pt} 
    \begin{tabular}{ll}
    \toprule
    \multicolumn{1}{p{13.25em}}{\textbf{Demographic Groups}} & \multicolumn{1}{p{7.915em}}{\textbf{Exposure rates}} \\
       \midrule
    \text{Female} & 0.058617 \\
    \midrule
    \text{Male} & 0.060259 \\
    \midrule
    \text{Age\_18\_29} & 0.052464 \\
    \midrule
    \text{Black\_non\_latinx} & 0.037505 \\
    \midrule
    \text{Age\_30\_39} & 0.067502 \\
    \midrule
    \text{Other\_race} & 0.097971 \\
    \midrule
    \text{Age\_40\_49} & 0.078672 \\
    \midrule
    \text{White\_non\_latinx} & 0.020698 \\
    \midrule
    \text{Age\_50\_59} & 0.077986 \\
    \midrule
    \text{Age\_60\_69} & 0.07898 \\
    \midrule
    \text{Age\_0\_17} & 0.038741 \\
    \midrule
    \text{Age\_70\_79} & 0.07275 \\
    \midrule
    \text{Two or more race} & 0.088107 \\
    \midrule
    \text{Age\_80\_} & 0.078768 \\
    \midrule
    \text{Asian\_non\_latinx} & 0.020833 \\
    \bottomrule
    \end{tabular}
    \end{adjustbox}
  \caption{\new{Exposure Rates of Population Subgroups-Chicago City}}\label{tab:rates}
\end{minipage}
\end{figure}

\new{Next, we perform a PCA analysis \cite{principal} using \attrib{Pop.} merged with the \attrib{Chicago-COVID-Zipcode} and produce a \emph{Biplot}, as shown in Figure~\ref{fig:hist}, to identify and compare the contribution of different attributes to the COVID-19 death and case numbers.} Note that we used cross-validation to tune the number of components of the PCA method.
To better visualize the impact of these attributes on different regions, we implement a $K$-means clustering \cite{Kmeans} method to partition our geographical areas (Zipcodes) based on the numbers of COVID-19 deaths and cases. We cluster them into three categories and label them as high, medium, and low impacted areas, accordingly.

Our findings in this part show that areas with higher COVID-19 deaths and cases tend to have a higher population, higher elderly population (compared with the young one), more Black Or African American, Latinos population (compared with other races), and lower median income. Furthermore, we believe that the revealed negative correlation between income and COVID-19 cases should raise serious concerns, in future decision-making procedures. One justification for that is lower-income individuals are mostly daily-paid and cannot afford living expenses if they are self-quarantined or stop working. Consequently, the pandemic inevitably affects the areas harder with larger lower-income populations.
However, we do not have COVID-19 data based on income level. Hence, we do not consider income in this study.

In brief, we can observe notable differences in contribution to the COVID-19 positive cases and death tolls across specific demographic attributes (e.g, older age). This fact reveals the critical role of a \emph{Fair-Diverse} model, which considers not only the overall population size but also the unrepresentative (exposed) population for scarce resource allocation problem. %\hadis{In our study, we assume that the income level difference is negligible across different regions, so we do not consider that in the construction of population subgroups. Is not it better to say income is a hidden impact of demographic attribute?}

%\end{wrapfigure}

\subsection{Estimating the Distribution of Exposed Population}\label{sec:risk}

Estimating the marginal distribution of the number of high-risk individuals in each demographic group, i.e., $P(E=1|g_i)$, requires the data of positive-tested (infected) cases and death tolls.
In this article, we intend to propose a resource allocation plan for vaccines and treatments for the COVID-19 pandemic.   
However, we only have access to the count of individuals of each social group who were infected from COVID-19, namely $P(g_i|E=1)$. 
The probability $P(E=1|g_i)$ can be calculated from the Bayes formula 
\[
P(E=1|g_i)= \frac{P(g_i|E=1)P(E=1)}{P(g_i)}.
\]

We will use $P(E=1|g_i)$, in the \S\ref{sec:res1} to form the fairness constraints of the optimization problem, Equation~\ref{eq:fairness}, and calculate the exposed population in each zip-code. \new{Table~\ref{tab:rates} presents the exposure rates for different population subgroups in Chicago City (The exposure rates for other cases (cities) are provided in Appendix~\ref{tab:rates_NYC},~\ref{tab:rates_BT}.} %The calculated exposure rates across different population subgroups will be applied to all zipcodes, similarly.
%\nazanin{need to justify why the exposed populations are different}

\subsection{Fair-Diverse Allocation}\label{sec:res1}

To evaluate the proposed \emph{Fair-Diverse} model, we construct \emph{Diverse-only}, \emph{Fair-only}, and equalized importance or \emph{alpha=0.5} models, and compare the allocation solutions as well as the resulted fairness and diversity gaps among them. In our terminology, \emph{Diverse-only} corresponds to an allocation that is merely based on diversity constraint, Equation~\ref{eq:diversity}, and is not considering any other (e.g., fairness) measures. Similarly, \emph{Fair-only} corresponds to a model solely based on fairness constraint. Furthermore, \emph{alpha=0.5} refers to a model that has equalized weights on Fairness and Diversity constraints in the optimization setting. 
% \new{We shall modify $b$ in each problem since the total capacity or total available units of vaccines to be allocated changes from one city to another.} 

As mentioned in \S~\ref{sec:opt}, we propose a diversity and fairness trade-off problem as in $P2$. As long as the resource constraint is satisfied, $P2$ has an optimal solution given an $\alpha$ value. However, the optimal solution obtained from $p2$ might not be feasible given the diversity and fairness requirements $\epsilon_d$ and $\epsilon_f$ as mentioned in \S~\ref{sec:feasibility}. That is we apply our proposed binary search algorithm to discover a range for $\alpha$ that results in a feasible solution.

% \nazanin{need to move some part of this discussion to the related work and put it before descriptive analysis}

To evaluate the performance of the Algorithm~\ref{alg:tune},\new{ we will consider three US cities and different sensitive attributes ( Race, Age, and Gender) and solve each \emph{instance} problem separately.} We will then show the allocation results for each problem using the four above-mentioned models and discuss $\alpha$ ranges in detail. \new{The total number of zip-code areas varies for each city with 177 regions in New York City, 58 in Chicago, and 36 in Baltimore}. The resulted allocations of each problem is sorted by population size for the top 15 area (zip-codes). 

\stitle{City of Chicago}

\new{We consider $b=200000$ units of vaccines as the total available resources to be allocated for the city of Chicago. Using \attrib{Pop.}, and \attrib{Chicago-COVID-Zipcode} datasets, we evaluate the allocation results of our proposed resource allocation framework at the zip-code level in the city of Chicago.}

The first instance problem considers \emph{Race} as the sensitive attribute and the associated \emph{Fair-Diverse} model captures inequalities across different \emph{racial subgroups}.
The results for the top 15 populated areas are reported in Table~\ref{tab:racial}. The baseline values for $\epsilon_{f}$ and $\epsilon_{d}$ are derived from the \emph{alpha=0.5} model and are equal to 0.24 and 0.007. Note that we do not use the tuning algorithm for this \emph{alpha=0.5} model. 
For \emph{Fair-Diverse} model $\epsilon_{f}$ and $\epsilon_{d}$ are both set to be 0.025 to decrease the fairness gap compared to the baseline value. 
The resulted tuned range for $\alpha$ is between 0.54 and 0.86 in this case (midpoint=0.70). 
Looking at Table~\ref{tab:racial}, the area associated with zip-code "$60639$" for instance, receives a lower number of vaccines using \emph{Diverse-only} model but higher in both \emph{alpha=0.5} and \emph{Fair-Diverse} model ($\alpha=0.70$) due to having higher total exposed population. In contrast, the \emph{Fair-only} closes the fairness gap to the full extent, ($\epsilon_{f}=0$), and as a result obtains an extreme allocation solution in which only a few areas (zip-codes) receive vaccines. Undoubtedly, this could not be a desirable allocation solution under certain fairness and diversity requirements. Since the table represents the top 15 populated areas, we cannot observe all areas with positive allocation using \emph{Fair-only} model.

\begin{table}[htbp]
  \centering
  \begin{adjustbox}{width=0.8\columnwidth}
    \begin{tabular}{rrrrrrr}
    \multicolumn{1}{l}{\textbf{Zipcode}} & \multicolumn{1}{l}{\textbf{Total population}} & \multicolumn{1}{l}{\textbf{Exposed Population}} & \multicolumn{1}{l}{\textbf{alpha=0.5}} & \multicolumn{1}{l}{\textbf{alpha tuned (0.70)}} & \multicolumn{1}{l}{\textbf{Diverse-only}} & \multicolumn{1}{l}{\textbf{Fair-only}} \\
    \midrule
    60629 & 113046 & 5802  & 10329 & 12125 & 9497  & 0 \\
    \midrule
    60618 & 91351 & 3652  & 7002  & 9798  & 7675  & 0 \\
    \midrule
    60623 & 91159 & 4815  & 8330  & 9777  & 7659  & 0 \\
    \midrule
    60639 & 89452 & 5145  & 8174  & 9594  & 7515  & 70901 \\
    \midrule
    60647 & 86586 & 3824  & 7912  & 9287  & 7274  & 0 \\
    \midrule
    60617 & 83590 & 3782  & 7638  & 8966  & 7023  & 0 \\
    \midrule
    60608 & 81930 & 3909  & 6280  & 8788  & 6883  & 0 \\
    \midrule
    60625 & 78085 & 2945  & 5986  & 8031  & 6560  & 0 \\
    \midrule
    60634 & 73894 & 2491  & 5664  & 4491  & 6208  & 0 \\
    \midrule
    60620 & 72094 & 2744  & 6587  & 7733  & 6057  & 0 \\
    \midrule
    60641 & 70992 & 2973  & 6455  & 7614  & 5964  & 0 \\
    \midrule
    60614 & 66485 & 1602  & 5096  & 4040  & 5586  & 0 \\
    \midrule
    60657 & 65841 & 1591  & 5047  & 4001  & 5532  & 0 \\
    \midrule
    60640 & 65412 & 2047  & 5014  & 3975  & 5496  & 0 \\
    \midrule
    60609 & 64420 & 3157  & 5886  & 6909  & 5412  & 0 \\
    \end{tabular}%
    \end{adjustbox}
    \vspace{2mm}
    \caption{{Resource Allocation (Racial groups): \new{Top 15 populated areas- Chicago}}}
  \label{tab:racial}%
\end{table}%

The second instance problem considers \emph{Age} as the sensitive attribute and the associated \emph{Fair-Diverse} model captures inequalities across different \emph{age groups}. The results for the top 15 populated areas are reported in Table~\ref{tab:age}. The baseline values for $\epsilon_{f}$ and $\epsilon_{d}$ are derived from the \emph{alpha=0.5} problem and are equal to 0.012 and 0. %The associated $\alpha$ range, in this case, is between zero and one. 
Next, we run the binary search algorithm to tune the $\alpha$ value and find a feasible solution for \emph{Fair-Diverse} model. In this case, $\epsilon_{f}$ and $\epsilon_{d}$ are both set to be 0.003. The resulted tuned range for $\alpha$ is between 0.57 and 1 (midpoint=0.78). As mentioned previously, the \emph{Diverse-only} model assigns vaccines to areas only based on the total population, and the \emph{Fair-only} model obtains an extreme allocation solution in which only a few areas (zip-codes) receive vaccines. Therefore, none of these models are capable of delivering a fair and diverse allocation solution. Note that in this instance problem, \emph{alpha=0.5} model is not doing any better than the \emph{Diverse-only} model since the weight on the fairness component is not adequate to change the results (diversity-gap is dominant). This result can further reveal the necessity of the proposed tuning algorithm.

\begin{table}[htbp]
  \centering
  \begin{adjustbox}{width=0.8\columnwidth}
    \begin{tabular}{rrrrrrr}
    \multicolumn{1}{l}{\textbf{Zipcode}} & \multicolumn{1}{l}{\textbf{Total population}} & \multicolumn{1}{l}{\textbf{Exposed Population}} & \multicolumn{1}{l}{\textbf{alpha=0.5}} & \multicolumn{1}{l}{\textbf{alpha tuned (0.78)}} & \multicolumn{1}{l}{\textbf{Diverse-only}} & \multicolumn{1}{l}{\textbf{Fair-only}} \\
   \midrule
    60629 & 113046 & 6367  & 9204  & 8989  & 9204  & 0 \\
    \midrule
    60618 & 91351 & 5478  & 7672  & 7492  & 7672  & 0 \\
    \midrule
    60623 & 91159 & 5084  & 7422  & 7249  & 7422  & 70300 \\
    \midrule
    60639 & 89452 & 5127  & 7367  & 7195  & 7367  & 0 \\
    \midrule
    60647 & 86586 & 5121  & 7312  & 7141  & 7312  & 0 \\
    \midrule
    60617 & 83590 & 5007  & 6931  & 7093  & 6931  & 0 \\
    \midrule
    60608 & 81930 & 4872  & 6904  & 6743  & 6904  & 0 \\
    \midrule
    60625 & 78085 & 4698  & 6569  & 6415  & 6569  & 0 \\
    \midrule
    60634 & 73894 & 4610  & 6219  & 6365  & 6219  & 0 \\
    \midrule
    60620 & 72094 & 4404  & 6009  & 6150  & 6009  & 0 \\
    \midrule
    60641 & 70992 & 4302  & 5955  & 5816  & 5955  & 0 \\
    \midrule
    60614 & 66485 & 4028  & 5735  & 5604  & 5735  & 0 \\
    \midrule
    60657 & 65841 & 4034  & 5730  & 5864  & 5730  & 0 \\
    \midrule
    60640 & 65412 & 4200  & 5644  & 5776  & 5644  & 0 \\
    \midrule
    60609 & 64420 & 3651  & 5264  & 5141  & 5264  & 0 \\
    \end{tabular}%
    \end{adjustbox}
    \vspace{2mm}
    \caption{{Resource Allocation (Age groups): \new{Top 15 populated areas-Chicago} }}
  \label{tab:age}%
\end{table}%

Finally, a noteworthy instance occurs when we consider \emph{gender} as the sensitive attribute and the \emph{Fair-Diverse} model attempts to eliminate the unfairness between male and female subgroups. The results for the top 15 populated areas are reported in Table~\ref{tab:gender}. The baseline values for $\epsilon_{f}$ and $\epsilon_{d}$ derived from the \emph{alpha=0.5} problem and are both close to zero ( 6.04e-05 and 0). It is worth mentioning that this is because of the similar gender population distribution across different regions. Consequently, the fairness and diversity requirements can be satisfied even with \emph{Diverse-only} model. We can still run the binary search algorithm to tune the $\alpha$ value and find a feasible solution for \emph{Fair-Diverse} model by closing the fairness gap further. To do this, the $\epsilon_{f}$ and $\epsilon_{d}$ values in \emph{Fair-Diverse} model should be set to 0 and 0.1. The resulted tuned range or $\alpha$ is between 0.92 and 1 in this case (midpoint=0.96). Looking at Table~\ref{tab:gender} and exposure rates (shown in Table~ \ref{tab:rates}), we notice that, in this specific instance problem, the exposed population size is, in actuality, aligned with the total population size in different areas. In other words, highly populated areas tend to have higher exposed populations as well. Therefore, the resulted allocations from \emph{Diverse-only}, \emph{alpha=0.5}, and \emph{Fair-Diverse} models are very close and even equal in some cases. 

\begin{table}[htbp]
  \centering
  \begin{adjustbox}{width=0.8\columnwidth}
    \begin{tabular}{rrrrrrr}
    \multicolumn{1}{l}{\textbf{Zipcode}} & \multicolumn{1}{l}{\textbf{Total population}} & \multicolumn{1}{l}{\textbf{Exposed Population}} & \multicolumn{1}{l}{\textbf{alpha=0.5}} & \multicolumn{1}{l}{\textbf{alpha tuned (0.96)}} & \multicolumn{1}{l}{\textbf{Diverse-only}} & \multicolumn{1}{l}{\textbf{Fair-only}} \\
   \midrule
    60629 & 113046 & 5802 & 9520  & 9529  & 9520  & 0 \\
    \midrule
    60618 & 91351 & 3651 & 7695  & 7703  & 7695  & 0 \\
    \midrule
    60623 & 91159 & 4815 & 7697  & 7705  & 7697  & 0 \\
    \midrule
    60639 & 89452 & 5145 & 7555  & 7563  & 7555  & 0 \\
    \midrule
    60647 & 86586 & 3824 & 7295  & 7302  & 7295  & 0 \\
    \midrule
    60617 & 83590 & 3782 & 7033  & 7026  & 7033  & 0 \\
    \midrule
    60608 & 81930 & 3909 & 6914  & 6921  & 6914  & 79366 \\
    \midrule
    60625 & 78085 & 2944 & 6573  & 6579  & 6573  & 0 \\
    \midrule
    60634 & 73894 & 2490 & 6209  & 6203  & 6209  & 0 \\
    \midrule
    60620 & 72094 & 2743 & 6035  & 6029  & 6035  & 0 \\
    \midrule
    60641 & 70992 & 2972 & 5989  & 5995  & 5989  & 0 \\
    \midrule
    60614 & 66485 & 1602 & 5567  & 5562  & 5567  & 0 \\
    \midrule
    60657 & 65841 & 1591 & 5515  & 5521  & 5515  & 0 \\
    \midrule
    60640 & 65412 & 2047 & 5498  & 5504  & 5498  & 0 \\
    \midrule
    60609 & 64420 & 3156 & 5424  & 5430  & 5424  & 0 \\
    \end{tabular}%
    \end{adjustbox}
    \vspace{2mm}
    \caption{{Resource Allocation (Gender groups): \new{Top 15 populated areas-Chicago}}}
  \label{tab:gender}%
\end{table}%

\begin{figure}[ht]
  \centering
  \begin{minipage}[b]{0.45\textwidth}
    \includegraphics[width=\textwidth]{{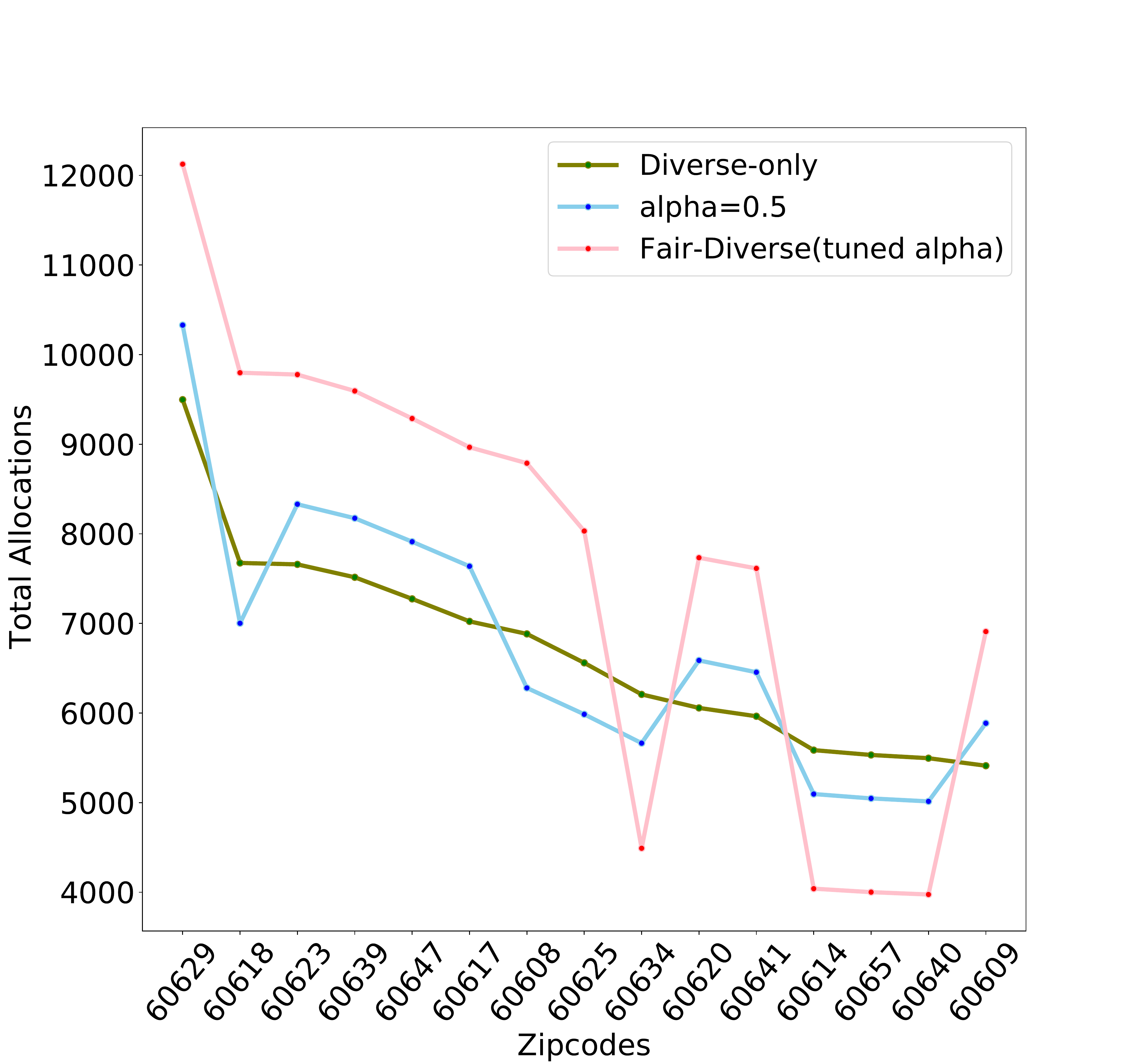}}
    \caption{Resource Allocation (Racial groups): \\ \new{Top 15 populated areas-Chicago}}
    \label{fig:subim1}
  \end{minipage}
  %\hfill
  \begin{minipage}[b]{0.45\textwidth}
    \includegraphics[width=\textwidth]{{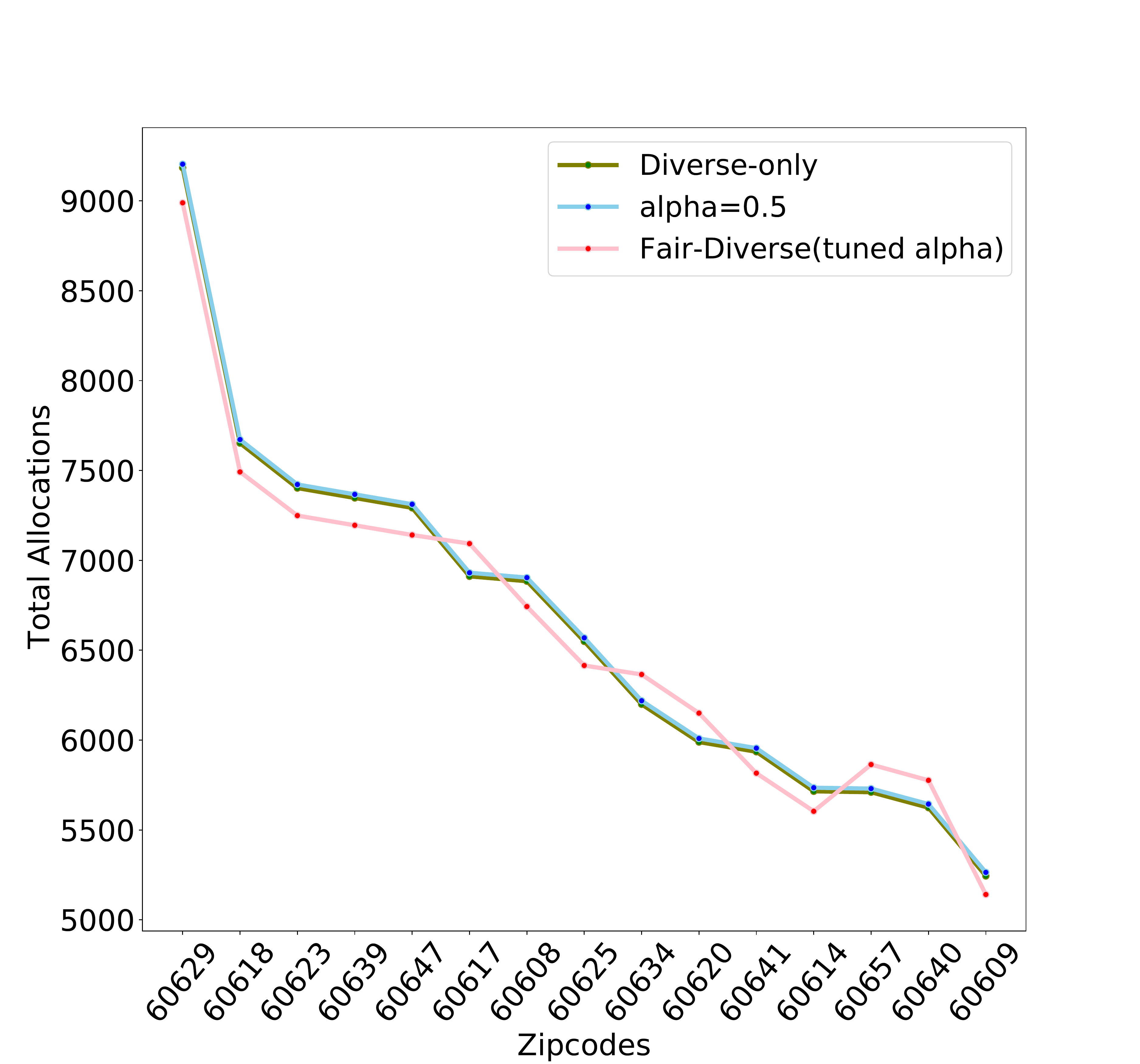}}
    \caption{Resource Allocation (Age groups): \\ \new{Top 15 populated areas-Chicago}}
        \label{fig:subim2}
  \end{minipage}
\end{figure}

Figures~\ref{fig:subim1} and \ref{fig:subim2} present the results for the top 15 populated regions in the Chicago City for \attrib{racial} and \attrib{age} instance problems. Note that in both instance problems, the total population equals the summation of all associated groups (e.g Age groups) due to having unknown labels in the data.
In Figure~\ref{fig:subim1}, for example, "60614" and "60634" regions receive less vaccines both under \emph{Fair-Diverse} and \emph{alpha=0.5} models due to having higher exposed population, Table~\ref{tab:racial}. The \emph{Diverse-only} does not consider the exposed population, therefore, the associated allocations are higher for these regions. Besides, the allocation obtained from the tuned range of $\alpha$ is significantly different from the allocation obtained with the \emph{alpha=0.5} model since the latter does not satisfy the fairness requirement $\epsilon_f$. 
Moving to another instance problem, Figure~\ref{fig:subim2} represents the allocation results for the age instance problem. Based on the plot, we can notice that the \emph{alpha=0.5} and \emph{Diverse-only} allocation solutions overlap. 
%However, we can observe different allocations for the  Fair-Diverse model with tuned $\alpha$ value. 
This can be justified by the fact that the 50\% emphasis on fairness is not sufficient to close the age subgroups disparities, and it requires a higher $\alpha$ value as we obtained through the tuning algorithm. That being said, the tuned $\alpha$ value, in this case, is 0.78, which is substantially higher than $0.5$.
For example, "60614" and "60609" regions receive less vaccines using \emph{Fair-Diverse} model with $\alpha$ tuning since they have relatively lower exposed population, Table~\ref{tab:age}. For observing more interactive visualization tools, please check our newly created web application on the Chicago City datasets using Rshiney\footnote{
\url{https://nazanin.shinyapps.io/Fair_Resource_Allocation/}}.

\new{\stitle{NewYork City (NYC)}}

\new{Utilizing \attrib{Pop.}, and \attrib{NYC-COVID-Zipcode} datasets, we assess our proposed resource allocation framework at zip-code level for New York City. In this problem, we consider the total number of vaccines available $b$ to 500000 since NYC has higher population compared with Chicago. Similar to the Chicago case study, the first instance problem for NYC considers \emph{Race}, and the second instance consider \emph{Age} as the sensitive attributes. We ignore the description of \emph{Gender} instance problem in this section (please see the appendix for the results).}

\new{The first instance problem considers \emph{Race} with the aim of capturing inequalities across different racial subgroups. The results for the top 15 populated areas in NYC are reported in Table~\ref{tab:NY_race}. The $\epsilon_{f}$ and $\epsilon_{d}$ derived from the \emph{alpha=0.5} model, are 0.10 and 0.0007 respectively. For \emph{Fair-Diverse} model $\epsilon_{f}$ and $\epsilon_{d}$ are both set to be 0.017 to decrease the fairness gap compared to the baseline value. The tuned $\alpha$ range is between 0.66 and 0.69 (midpoint=0.67). Looking at Table~\ref{tab:NY_race}, the area associated with zip-code "10467", receives a lower number of vaccines using \emph{Diverse-only} model but higher numbers in both \emph{alpha=0.5} and \emph{Fair-Diverse} model ($\alpha=0.67$) due to having higher total exposed population and more risks. In contrast, the \emph{Fair-only} closes the fairness gap to the full extent ($\epsilon_{f}=0$), and as a result, obtains an extreme allocation solution in which only a few areas (zip-codes) receive vaccines. As shown in Table~\ref{tab:NY_race}, none of the top 15 populated areas will receive vaccines under this extreme condition. As a result, this could not be an applicable allocation solution.} 

\new{The second instance problem considers \emph{Age} to captures inequalities across different \emph{age groups}. The results for the top 15 populated areas are reported in Table~\ref{tab:NY_age}. The baseline values for $\epsilon_{f}$ and $\epsilon_{d}$ are derived from the \emph{alpha=0.5} problem and are equal to 0.008 and 0 while the $\alpha$ value range is between zero and one. Next, we run the binary search algorithm to tune the $\alpha$ value and find a feasible solution for \emph{Fair-Diverse} model with $\epsilon_{f}$ and $\epsilon_{d}$ both set to 0.003. The resulted tuned range for $\alpha$ is between 0.68 and 1 (midpoint=0.84). As mentioned previously, the \emph{Diverse-only} model assigns vaccines to areas merely based on the total population, and the \emph{Fair-only} model obtains an extreme allocation solution in which only a few areas (zip-codes) receive vaccines. Therefore, these models are not capable of delivering a fair and diverse allocation solution. Note that in this instance problem, \emph{alpha=0.5} model is not doing any better than the \emph{Diverse-only} model since the weight on the fairness component is not adequate which further reveals the necessity of the proposed tuning approach.}

\begin{figure}[ht]
  \centering
  \begin{minipage}[b]{0.45\textwidth}
    \includegraphics[width=\textwidth]{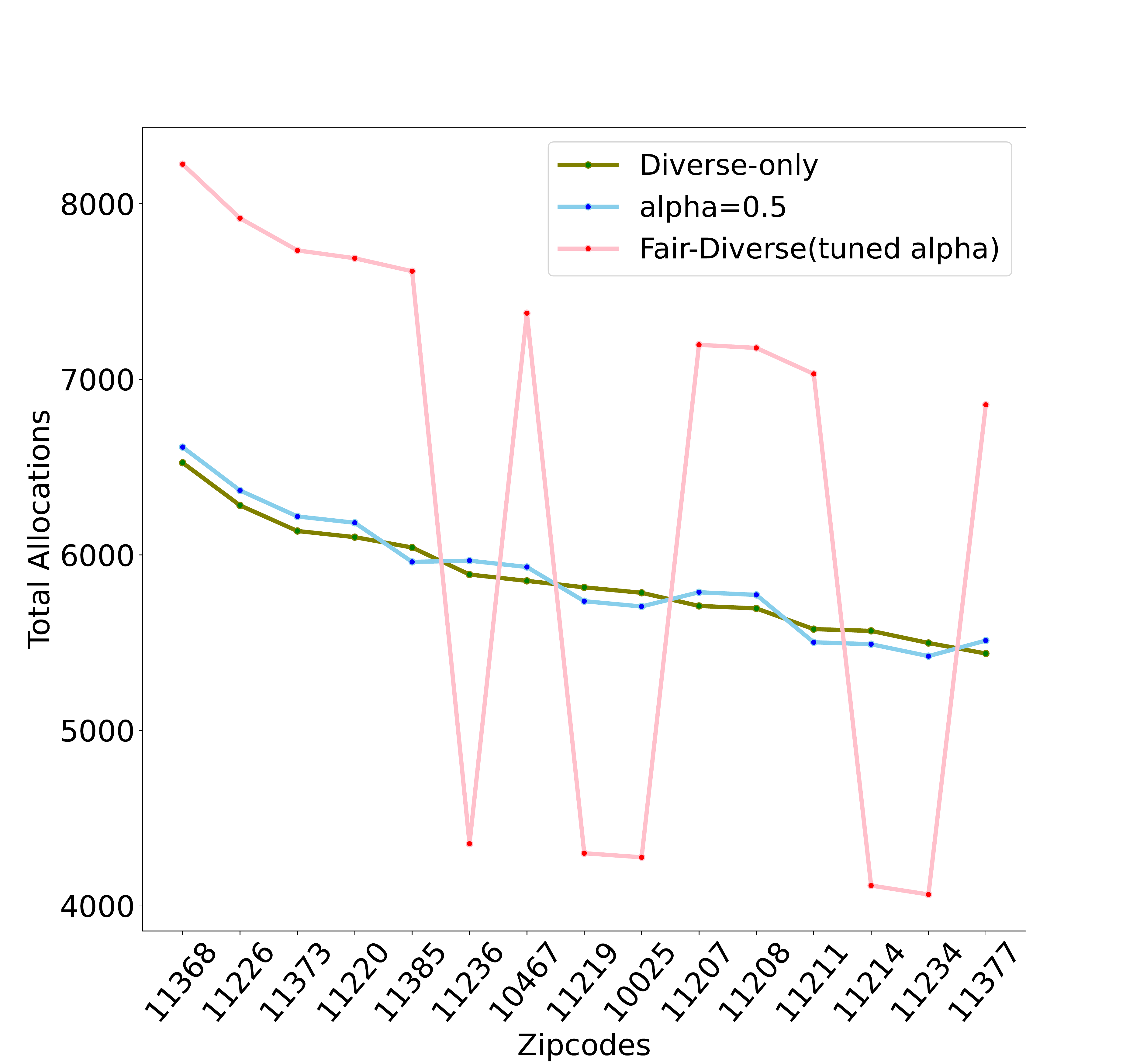}
    \caption{Resource Allocation (Racial groups): \\ \new{Top 15 populated areas-NYC}}
    \label{fig:NY_race}
  \end{minipage}
  %\hfill
  \begin{minipage}[b]{0.45\textwidth}
    \includegraphics[width=\textwidth]{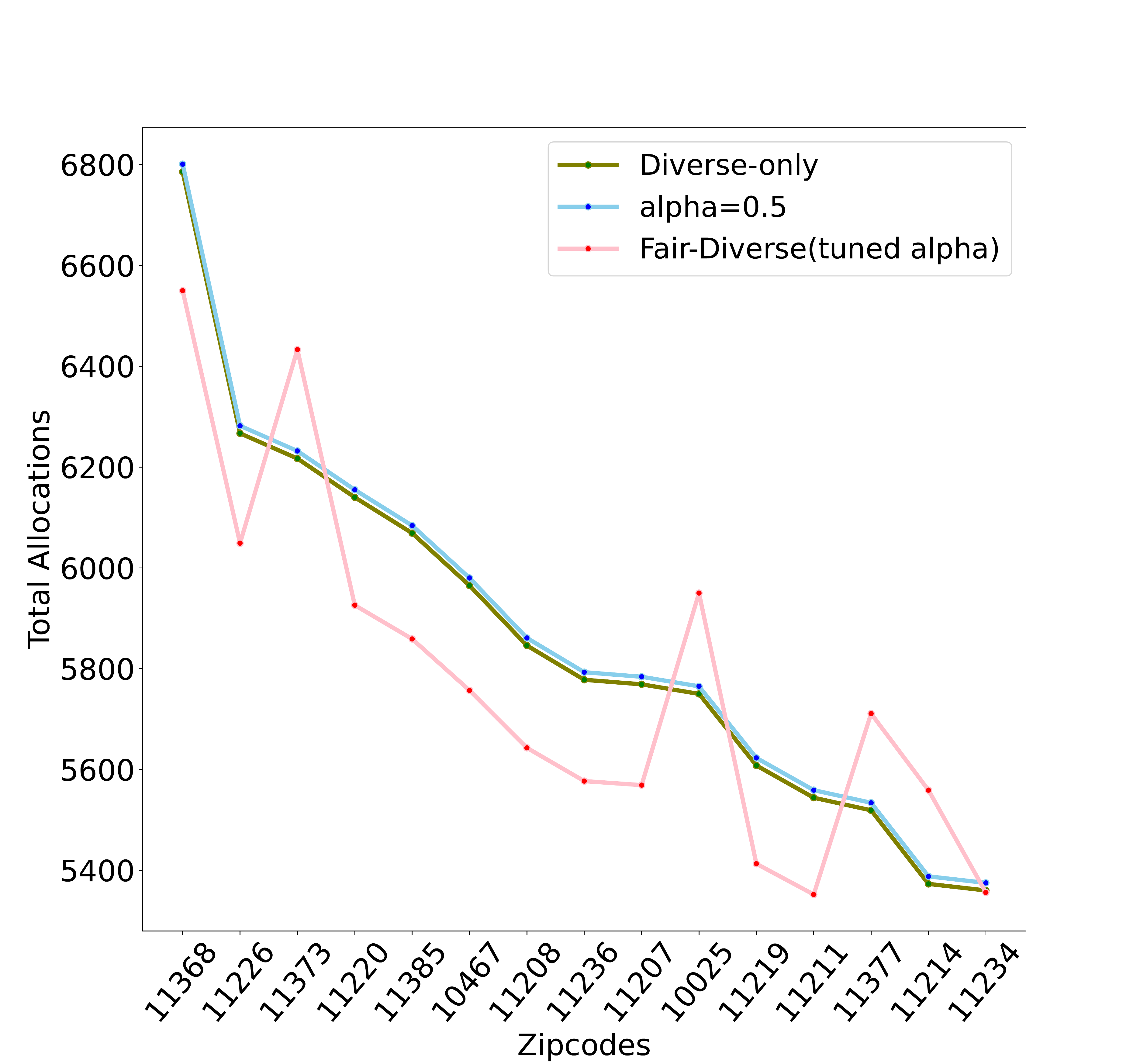}
    \caption{Resource Allocation (Age groups): \\ \new{Top 15 populated areas-NYC}}
        \label{fig:NY_age}
  \end{minipage}
\end{figure}

\begin{table}[htbp]
  \centering
\begin{adjustbox}{width=0.8\columnwidth}
    \begin{tabular}{ccccccc}
    \textbf{Zipcode} & \textbf{Total Population} & \textbf{Exposed Population} & \textbf{alpha=0.5} & \textbf{alpha\_tuned} & \textbf{Diverse-only} & \textbf{Fair-Only} \\
    \midrule
    11368 & 101558 & 12156 & 6614  & 8226  & 6525  & 0 \\
    \midrule
    11226 & 97766 & 7186  & 6367  & 7918  & 6282  & 0 \\
    \midrule
    11373 & 95500 & 8765  & 6219  & 7735  & 6136  & 0 \\
    \midrule
    11220 & 94949 & 8798  & 6183  & 7690  & 6101  & 0 \\
    \midrule
    11385 & 94032 & 7941  & 5960  & 7616  & 6042  & 0 \\
    \midrule
    11236 & 91632 & 6051  & 5967  & 4354  & 5888  & 0 \\
    \midrule
    10467 & 91077 & 8504  & 5931  & 7377  & 5852  & 0 \\
    \midrule
    11219 & 90499 & 5798  & 5736  & 4300  & 5815  & 0 \\
    \midrule
    10025 & 90025 & 5950  & 5706  & 4277  & 5784  & 0 \\
    \midrule
    11207 & 88857 & 7579  & 5787  & 7197  & 5709  & 0 \\
    \midrule
    11208 & 88639 & 8445  & 5772  & 7179  & 5695  & 0 \\
    \midrule
    11211 & 86805 & 6049  & 5502  & 7031  & 5577  & 0 \\
    \midrule
    11214 & 86639 & 5843  & 5491  & 4116  & 5567  & 0 \\
    \midrule
    11234 & 85567 & 5208  & 5423  & 4065  & 5498  & 0 \\
    \midrule
    11377 & 84635 & 7242  & 5512  & 6855  & 5438  & 0 \\
    \end{tabular}%
    \end{adjustbox}
     \caption{{Resource Allocation (Racial groups): \new{Top 15 populated areas- NYC}}}
  \label{tab:NY_race}%
\end{table}%

\new{Figures~\ref{fig:NY_race} and \ref{fig:NY_age} present the results for the top 15 populated regions in the New York City for \attrib{racial} and \attrib{age} instance problems. Note that in both instance problems, the total population equals the summation of all associated groups (e.g Age groups) due to having unknown labels in the data.
In Figure~\ref{fig:NY_race}, "11236" and "11219" zip-codes receive less vaccines both under \emph{Fair-Diverse} and \emph{alpha=0.5} models due to having lower exposed population, Table~\ref{tab:NY_race}. The \emph{Diverse-only} does not consider the exposed population, therefore, the associated allocations are higher for these regions. Besides, the allocation obtained from the tuned range of $\alpha$ is significantly different from the allocation obtained with the \emph{alpha=0.5} model since the latter does not satisfy the fairness requirement $\epsilon_f$.} 

\new{Moving to another instance problem, Figure~\ref{fig:NY_age} represents the allocation results for the age instance problem. Based on the plot, we can notice that the \emph{alpha=0.5} and \emph{Diverse-only} allocation solutions overlap. This can be justified by the fact that the 50\% emphasis on fairness is not sufficient to close the age subgroups disparities, and it requires a higher $\alpha$ value as we obtained through the tuning algorithm. That being said, the tuned $\alpha$ value, in this case, is 0.84, which is substantially higher than $0.5$.
For example, "11226" and "11211" regions receive less vaccines using \emph{Fair-Diverse} model comparing with \emph{Diverse-only} or \emph{alpha=0.5} since they have relatively lower exposed population, Table~\ref{tab:NY_age}.}

\begin{table}[htbp]
  \centering
\begin{adjustbox}{width=0.8\columnwidth}
    \begin{tabular}{ccccccc}
    \textbf{Zipcode} & \textbf{Total Population} & \textbf{Exposed Population} & \textbf{alpha=0.5} & \textbf{alpha\_tuned} & \textbf{Diverse-only} & \textbf{Fair-Only} \\
    \midrule
    11368 & 108904 & 10192   & 6786  & 6550  & 6786  & 0 \\
    \midrule
    11226 & 100574 & 9711  & 6267  & 6049  & 6267  & 0 \\
    \midrule
    11373 & 99772 & 9835  & 6217  & 6433  & 6217  & 0 \\
    \midrule
    11220 & 98531 & 9403  & 6140  & 5926  & 6140  & 0 \\
    \midrule
    11385 & 97405 & 9466  & 6069  & 5859  & 6069  & 0 \\
    \midrule
    10467 & 95723 & 9076  & 5965  & 5757  & 5965  & 0 \\
    \midrule
    11208 & 93812 & 8674  & 5846  & 5643  & 5846  & 0 \\
    \midrule
    11236 & 92721 & 9020 & 5778  & 5577  & 5778  & 0 \\
    \midrule
    11207 & 92591 & 8633  & 5769  & 5569  & 5769  & 0 \\
    \midrule
    10025 & 92284 & 9549  & 5750  & 5950  & 5750  & 0 \\
    \midrule
    11219 & 90004 & 8005  & 5608  & 5413  & 5608  & 75932 \\
    \midrule
    11211 & 88978 & 8286  & 5544  & 5352  & 5544  & 0 \\
    \midrule
    11377 & 88572 & 8808  & 5519  & 5711  & 5519  & 0 \\
    \midrule
    11214 & 86223 & 8725 & 5373  & 5559  & 5373  & 0 \\
    \midrule
    11234 & 86012 & 8536  & 5360  & 5356  & 5360  & 0 \\
    \end{tabular}%
    \end{adjustbox}
     \caption{{Resource Allocation (Age groups): \new{Top 15 populated areas- NYC}}}
  \label{tab:NY_age}%
\end{table}%

\new{\stitle{Baltimore City}}

\new{Utilizing \attrib{Pop.}, and \attrib{BLT-COVID-Zipcode} datasets, we empirically tested our proposed resource allocation framework at the zip-code level. In this problem, we modify the total number of vaccines available $b = 100000$ since Baltimore has a relatively lower population compared to Chicago and New York City. However, similar
to other case studies, the first instance problem considers \emph{Race}, and the second instance considers \emph{Age} as the sensitive attributes and ignore the description of Gender instance problem (please see the appendix for the results).}

\new{The \emph{Race} instance problem captures the inequalities across different racial subgroups. The results for the top 15 populated areas in Baltimore are reported in Table~\ref{tab:BT_race}. The $\epsilon_{f}$ and $\epsilon_{d}$ derived from the \emph{alpha=0.5} model, are 0.046 and 0.016 respectively. For \emph{Fair-Diverse} model $\epsilon_{f}$ and $\epsilon_{d}$ are both set to be 0.025 to decrease the fairness gap compared to the baseline value. The tuned $\alpha$ range is between 0.57 and 1 (midpoint=0.78). Looking at Table~\ref{tab:BT_race}, the area associated with zip-code "21230", receives a higher number of vaccines using \emph{Diverse-only} model but lower numbers in both \emph{alpha=0.5} and \emph{Fair-Diverse} model ($\alpha=0.78$) due to having lower total exposed population and less risk level. On the other hand, the \emph{Fair-only} closes the fairness gap ($\epsilon_{f}=0$), and obtains an extreme allocation solution in which only a few areas (zip-codes) receive vaccines. As shown in Table~\ref{tab:BT_race}, four regions among the top 15populated areas receive vaccines under this extreme condition. As a result, this could not be an applicable allocation solution.} 

\new{The \emph{Age} instance problem captures inequalities across different \emph{age groups}. The results for the top 15 populated areas are reported in Table~\ref{tab:BT_age}. The baseline values for $\epsilon_{f}$ and $\epsilon_{d}$ are derived from the \emph{alpha=0.5} problem and are equal to 0.014 and 0 while the $\alpha$ value range is between zero and one. Next, we run the binary search algorithm to tune the $\alpha$ value and find a feasible solution for \emph{Fair-Diverse} model with $\epsilon_{f}$ and $\epsilon_{d}$ both set to 0.007. The resulted tuned range for $\alpha$ is between 0.76 and 1 (midpoint=0.88). As mentioned previously, the \emph{Diverse-only} model assigns vaccines to areas merely based on the total population, and the \emph{Fair-only} model obtains an extreme allocation solution in which only four areas (zip-codes) receive vaccines. Therefore, these models are not capable of delivering a fair and diverse allocation solution. Note that in this instance problem, \emph{alpha=0.5} model is not doing any better than the \emph{Diverse-only} model since the weight on the fairness component is not adequate which further reveals the necessity of the proposed tuning approach.} 

\begin{figure}[ht]
  \centering
  \begin{minipage}[b]{0.45\textwidth}
    \includegraphics[width=\textwidth]{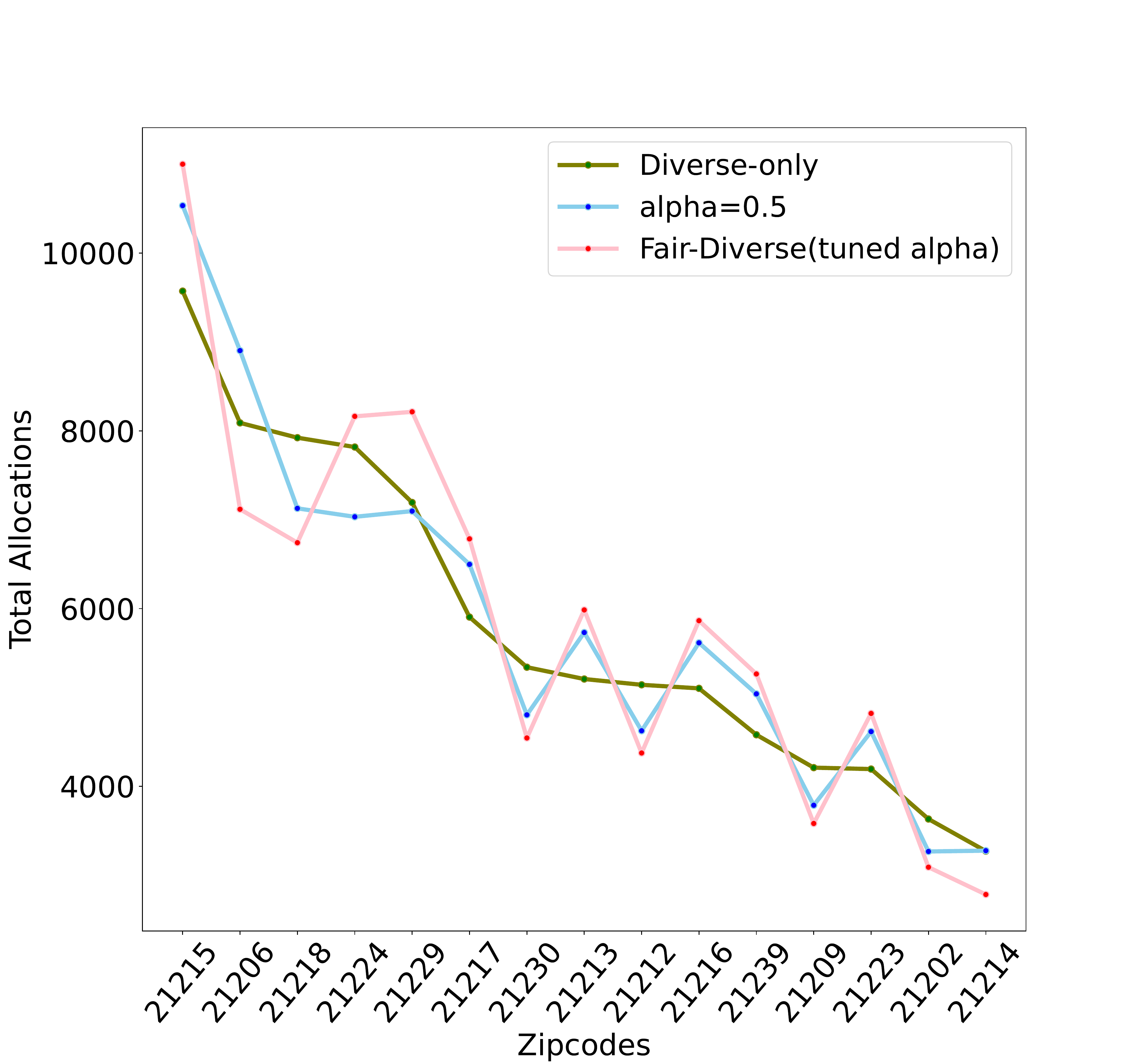}
    \caption{Resource Allocation (Racial groups): \\\new{top 15 populated areas-Baltimore}}
    \label{fig:BT_race}
  \end{minipage}
  %\hfill
  \begin{minipage}[b]{0.45\textwidth}
    \includegraphics[width=\textwidth]{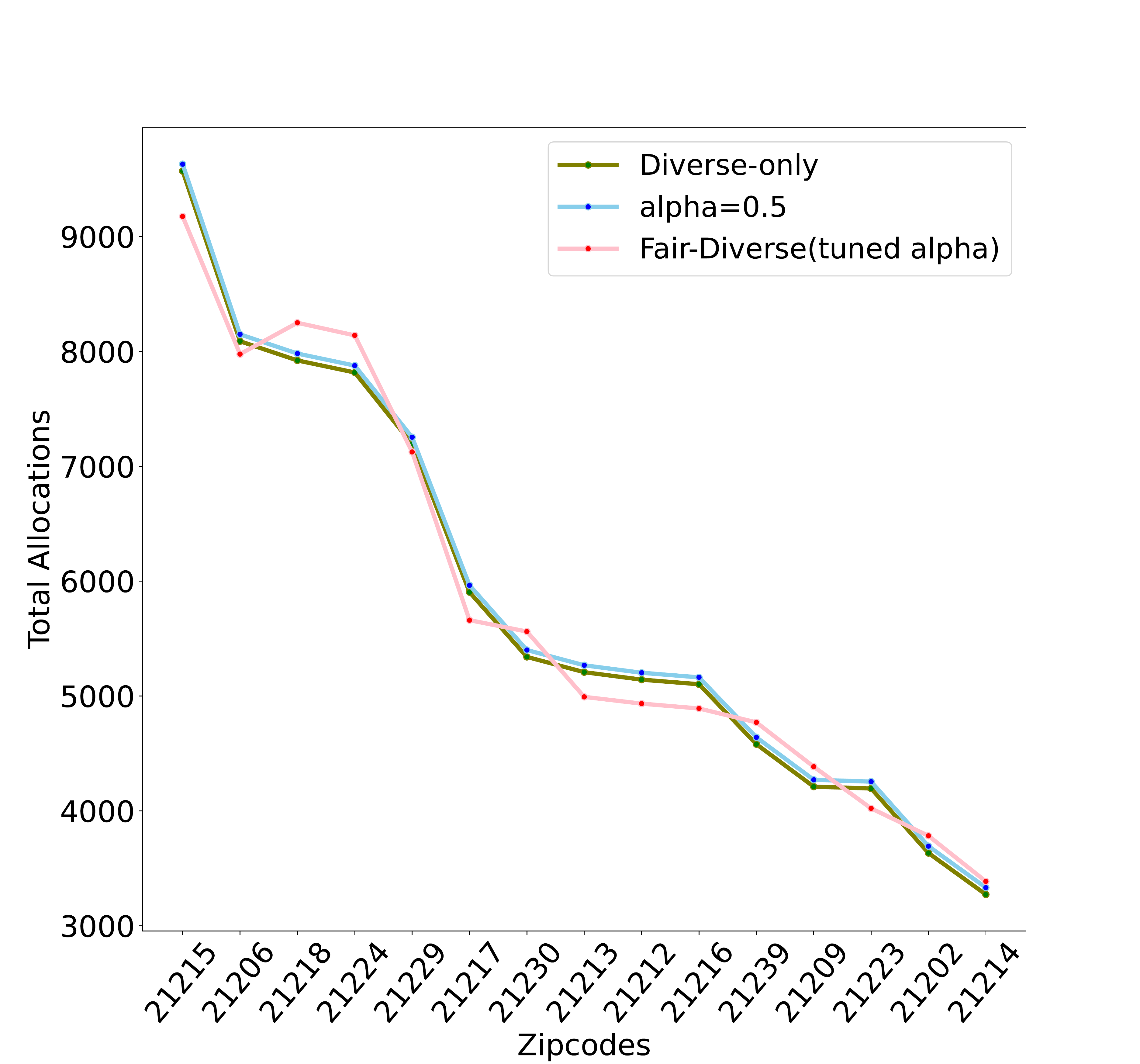}
    \caption{Resource Allocation (Age groups): \\ \new{top 15 populated areas-Baltimore}}
        \label{fig:BT_age}
  \end{minipage}
\end{figure}

% Table generated by Excel2LaTeX from sheet 'Baltimore-Race'
\begin{table}[htbp]
  \centering
    \begin{adjustbox}{width=0.8\columnwidth}
    \begin{tabular}{ccccccc}
    \textbf{Zipcode} & \textbf{Total Population} & \textbf{Exposed Population} & \textbf{alpha=0.5} & \textbf{alpha\_tuned} & \textbf{Diverse-only} & \textbf{Fair-Only} \\
    \midrule
    21215 & 60161 & 4893  & 10534 & 11000 & 9572  & 0 \\
    \midrule
    21206 & 50846 & 3944  & 8903  & 7118  & 8090  & 0 \\
    \midrule
    21218 & 49796 & 3739  & 7128  & 6742  & 7923  & 0 \\
    \midrule
    21224 & 49134 & 3782  & 7033  & 8163  & 7818  & 14868 \\
    \midrule
    21229 & 45213 & 3555  & 7097  & 8214  & 7194  & 12266 \\
    \midrule
    21217 & 37111 & 3050  & 6498  & 6785  & 5905  & 0 \\
    \midrule
    21230 & 33568 & 2188  & 4805  & 4545  & 5341  & 0 \\
    \midrule
    21213 & 32733 & 2722  & 5731  & 5985  & 5208  & 0 \\
    \midrule
    21212 & 32322 & 2176  & 4626  & 4376  & 5143  & 0 \\
    \midrule
    21216 & 32071 & 2725  & 5615  & 5864  & 5103  & 49078 \\
    \midrule
    21239 & 28793 & 2303  & 5041  & 5265  & 4581  & 0 \\
    \midrule
    21209 & 26465 & 1494  & 3788  & 3583  & 4211  & 11242 \\
    \midrule
    21223 & 26366 & 2143  & 4616  & 4821  & 4195  & 0 \\
    \midrule
    21202 & 22832 & 1724  & 3268  & 3091  & 3633  & 0 \\
    \midrule
    21214 & 20564 & 1510  & 3277  & 2784  & 3272  & 0 \\
    \end{tabular}%
    \end{adjustbox}
     \caption{{Resource Allocation (Racial groups): \new{Top 15 populated areas- Baltimore}}}
  \label{tab:BT_race}%
\end{table}%

\new{Figures~\ref{fig:BT_race} and \ref{fig:BT_age} show the results for the top 15 populated regions in the Baltimore City for \attrib{racial} and \attrib{age} instance problems. Note that in both instance problems, the total population equals the summation of all associated groups (e.g Age groups) due to having unknown labels in the data.
In Figure~\ref{fig:BT_race}, and zip-codes receive less vaccines both under \emph{Fair-Diverse} and  \emph{alpha=0.5} models due to having lower exposed population, Table~\ref{tab:BT_race}. The \emph{Diverse-only} does not consider the exposed population, therefore, the associated allocations are higher for these regions. Besides, the allocation obtained from the tuned range of $\alpha$ is significantly different from the allocation obtained with the \emph{alpha=0.5} model since the latter does not satisfy the fairness requirement $\epsilon_f$.} 

\new{Moving to the next instance problem, Figure~\ref{fig:BT_age} represents the allocation results for the age instance problem. Based on the plot, we can notice that the \emph{alpha=0.5} and \emph{Diverse-only} allocation solutions overlap. This can be justified by the fact that the 50\% emphasis on fairness is not sufficient to close the age subgroups disparities, and it requires a higher $\alpha$ value as we obtained through the tuning algorithm. That being said, the tuned $\alpha$ value, in this case, is 0.88, which is substantially higher than $0.5$.
For example, "21217" and "21223" regions receive less vaccines using \emph{Fair-Diverse} model comparing with \emph{Diverse-only} or \emph{alpha=0.5} since they have relatively lower exposed population, Table~\ref{tab:BT_age}.}

% Table generated by Excel2LaTeX from sheet 'Baltimore-Age'
\begin{table}[htbp]
  \centering
    \begin{adjustbox}{width=0.8\columnwidth}
    \begin{tabular}{ccccccc}
    \textbf{Zipcode} & \textbf{Total Population} & \textbf{Exposed Population} & \textbf{alpha=0.5} & \textbf{alpha\_tuned} & \textbf{Diverse-only} & \textbf{Fair-Only} \\
    \midrule
    21215 & 60161 & 4943  & 9572  & 9177  & 9572  & 13765 \\
    \midrule
    21206 & 50846 & 4211  & 8090  & 7978  & 8090  & 0 \\
    \midrule
    21218 & 49796 & 4213  & 7923  & 8251  & 7923  & 0 \\
    \midrule
    21224 & 49134 & 4285  & 7818  & 8141  & 7818  & 0 \\
    \midrule
    21229 & 45213 & 3762  & 7194  & 7126  & 7194  & 0 \\
    \midrule
    21217 & 37111 & 3040  & 5905  & 5661  & 5905  & 0 \\
    \midrule
    21230 & 33568 & 2961  & 5341  & 5562  & 5341  & 0 \\
    \midrule
    21213 & 32733 & 2661  & 5208  & 4993  & 5208  & 7954 \\
    \midrule
    21212 & 32322 & 2682  & 5143  & 4934  & 5143  & 0 \\
    \midrule
    21216 & 32071 & 2614  & 5103  & 4892  & 5103  & 0 \\
    \midrule
    21239 & 28793 & 2416  & 4581  & 4771  & 4581  & 0 \\
    \midrule
    21209 & 26465 & 2230  & 4211  & 4385  & 4211  & 2268 \\
    \midrule
    21223 & 26366 & 2150  & 4195  & 4022  & 4195  & 0 \\
    \midrule
    21202 & 22832 & 2048  & 3633  & 3783  & 3633  & 0 \\
    \midrule
    21214 & 20564 & 1727  & 3272  & 3386  & 3272  & 26707 \\
    \end{tabular}%
    \end{adjustbox}
     \caption{{Resource Allocation (Age groups): \new{Top 15 populated areas- Baltimore}}}
  \label{tab:BT_age}%
\end{table}%

As we discussed in \S~\ref{sec:feasibility}, the allocation solution that is obtained from $P2$ does not necessarily satisfy the fairness and diversity requirements (under any $\alpha$ value). To demonstrate the performance of the tuning algorithm, which always returns a range for $\alpha$ under which the optimal solution of $P2$ is feasible, we now study the impact of $\epsilon_f$ and $\epsilon_d$ on the optimal solution under different models using the city of Chicago case study.
%\new{For the purpose of this analysis, we present the results of this part based on the city of Chicago data for it is a diverse city with a relatively high Hispanic and Black (vulnerable) population}.  

\begin{figure}
\centering
\subfigure[]{\label{fig:5_1}\includegraphics[width=45mm]{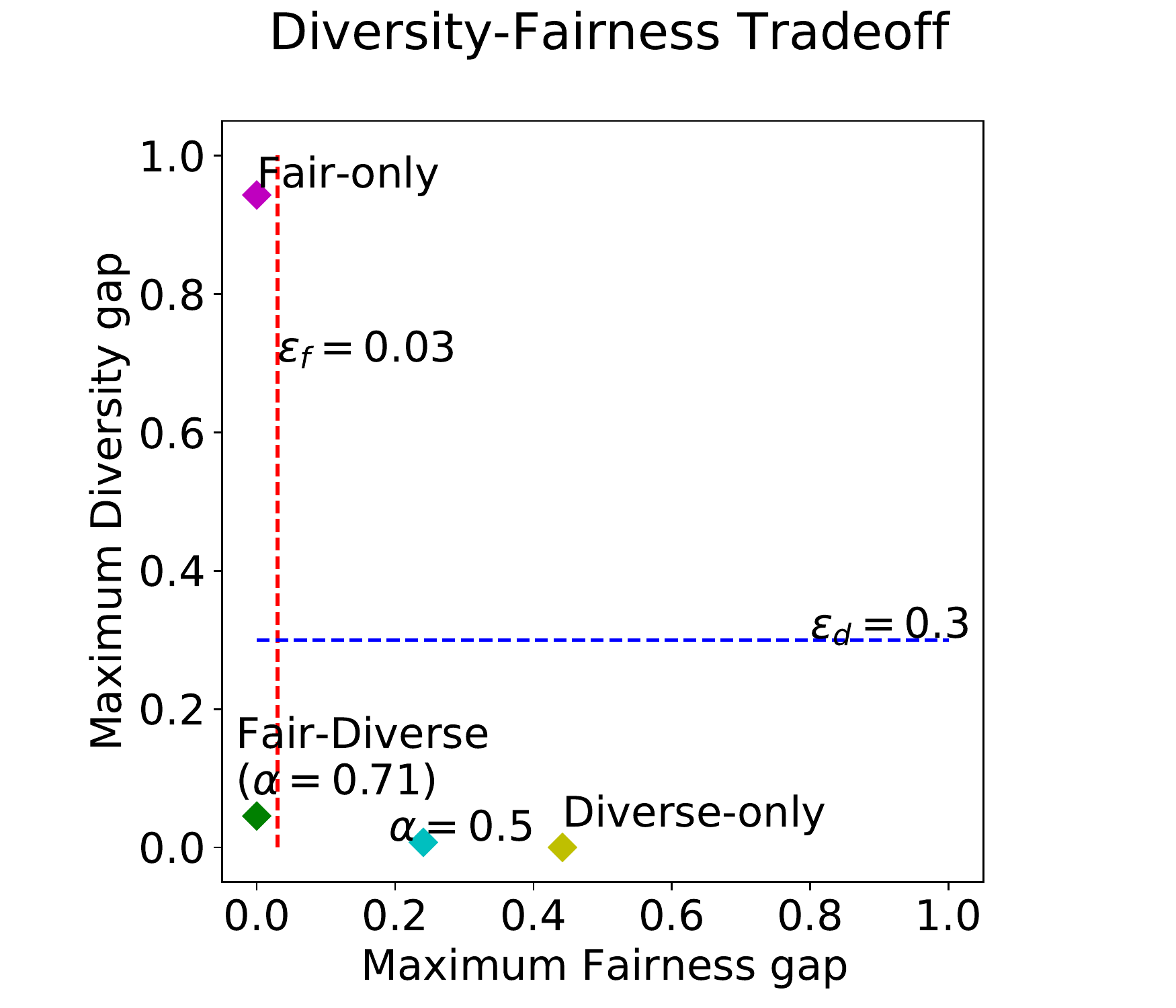}}
\subfigure[]{\label{fig:5_2}\includegraphics[width=45mm]{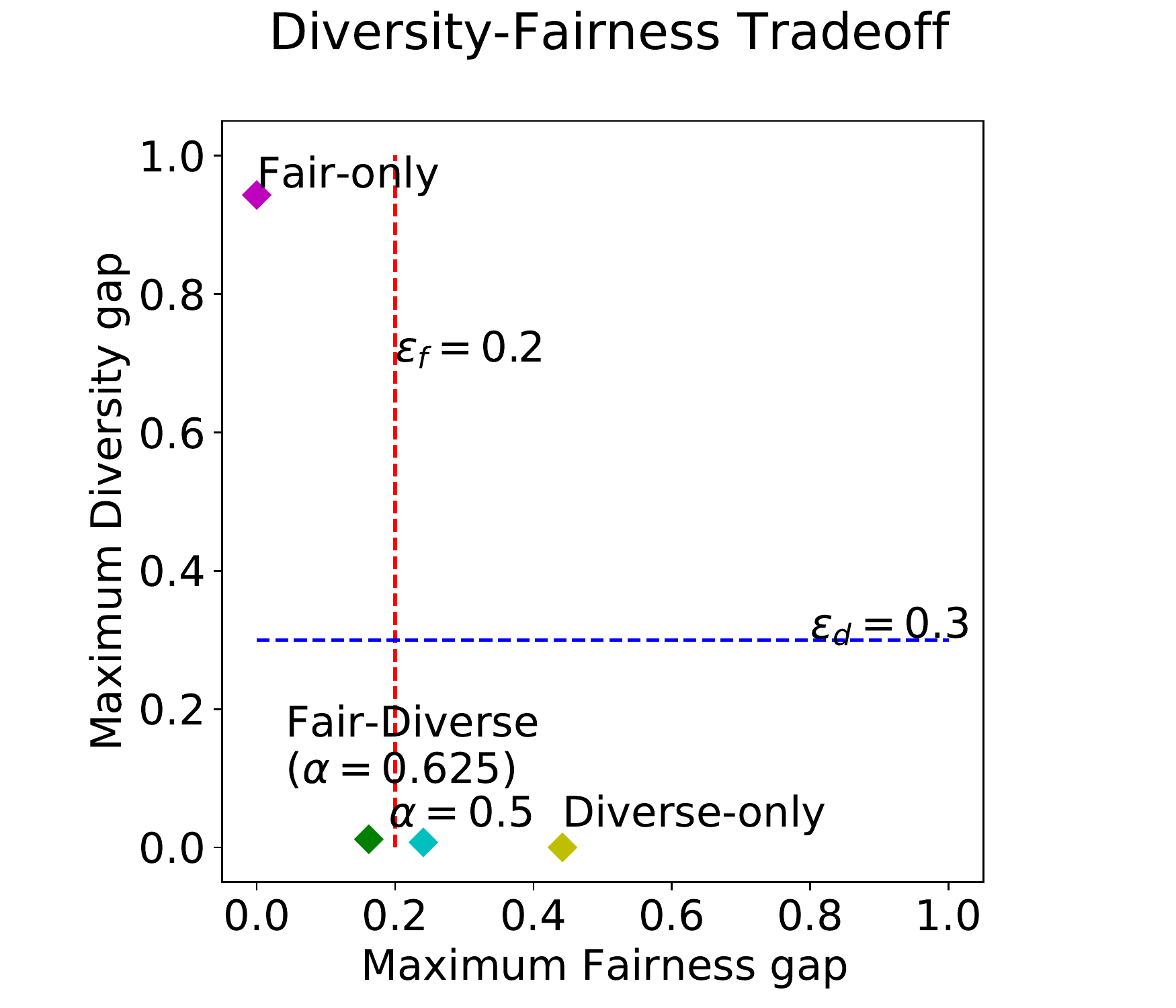}}
\subfigure[]{\label{fig:5_3}\includegraphics[width=45mm]{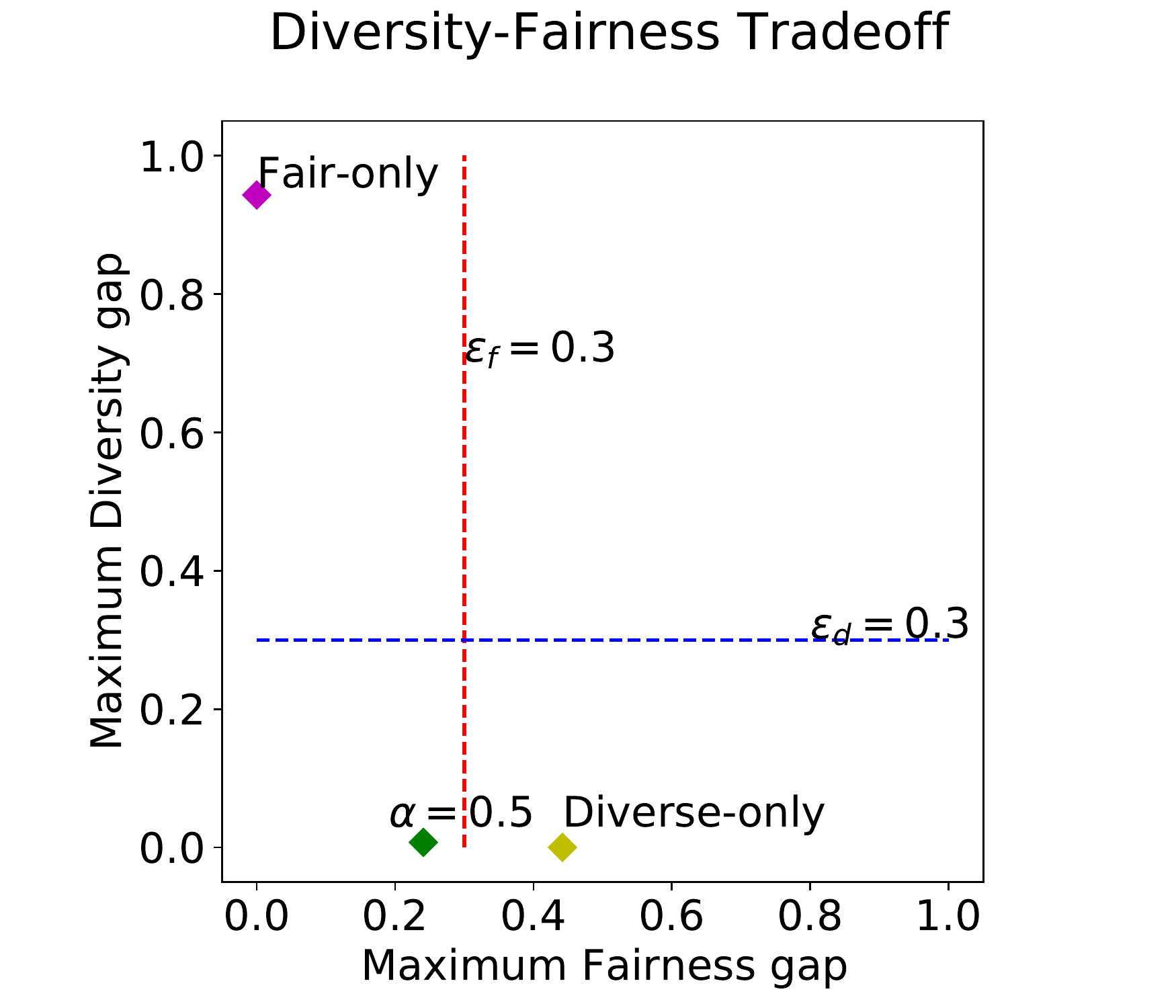}}
\caption{Impact of the $\epsilon_f$ and $\epsilon_d$ on the diversity-fairness trade-off}\label{fig:5}
\end{figure}

The plots in \ref{fig:5} are based on the \emph{racial} instance problem. These figures reveal that under any fairness and diversity requirements (given $\epsilon_f$ and $\epsilon_d$ values) the $\alpha$ tuning algorithm returns a feasible solution for $\emph{Fair-Diverse}$ model. Note that, this is not the case for other models (\emph{Fair-only}, \emph{Diverse-only} and \emph{alpha=0.5}) as it can be observed from the Figures~\ref{fig:5_1} and~\ref{fig:5_2}. In other words, the optimal solutions obtained from the \emph{Diverse-only} and \emph{alpha=0.5} models does not satisfy the fairness requirement ($\epsilon_{f} \leq 0.03$ and $\epsilon_{f} \leq 0.2$ ), and the solution obtained from the \emph{Fair-only} model fails to satisfy the diversity requirement ($\epsilon_{d} \leq 0.3$) in Figures~\ref{fig:5_1} and \ref{fig:5_2}. However, if we relax the fairness requirement to $\epsilon_f \leq 0.3$, Figure~\ref{fig:5_3}, the \emph{alpha=0.5} model can indeed achieve a feasible solution. It is worth mentioning that the diversity requirement is easier to achieve compared to the fairness requirement due to the larger inherent disparities in exposed population. %As a result, the

\subsection{Price of Fairness}
In this section, we will compare the fairness and diversity gaps under different models and population subgroups to discuss the price of fairness (PoF). \new{The results are obtained based on the aforementioned racial group instance problem for the city of Chicago. We present the results under tuned $\alpha$ value (0.71) in this part.}

\begin{figure}[!htbp]
    \begin{minipage}{0.3\textwidth}
        \includegraphics[width=\linewidth]{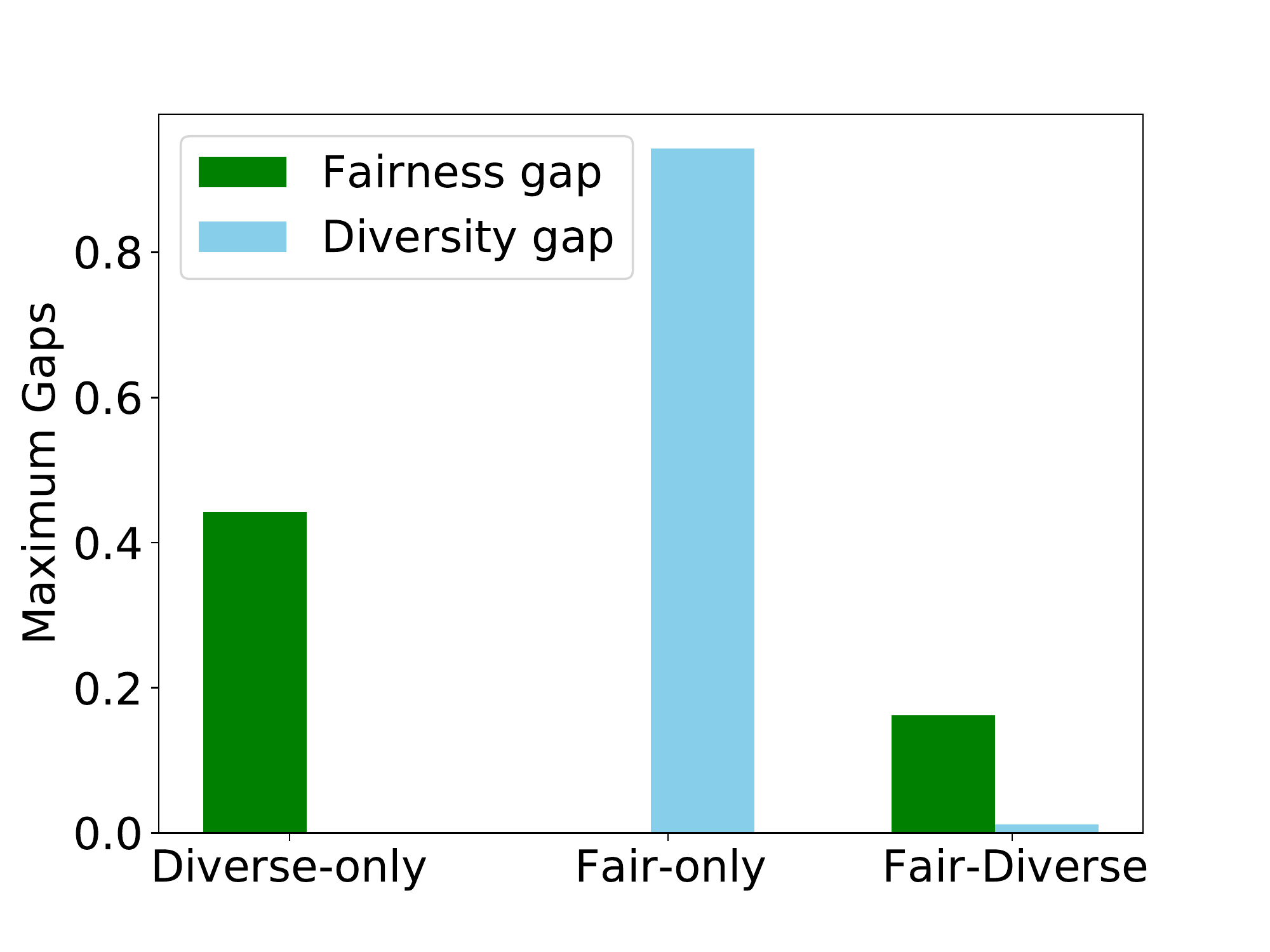}
    %\vspace{-5mm}
    \caption{\small Diversity and Fairness gaps in different models}
    \label{fig:gaps1}
    \end{minipage}
        \begin{minipage}{0.33\textwidth}
        \includegraphics[width=\linewidth]{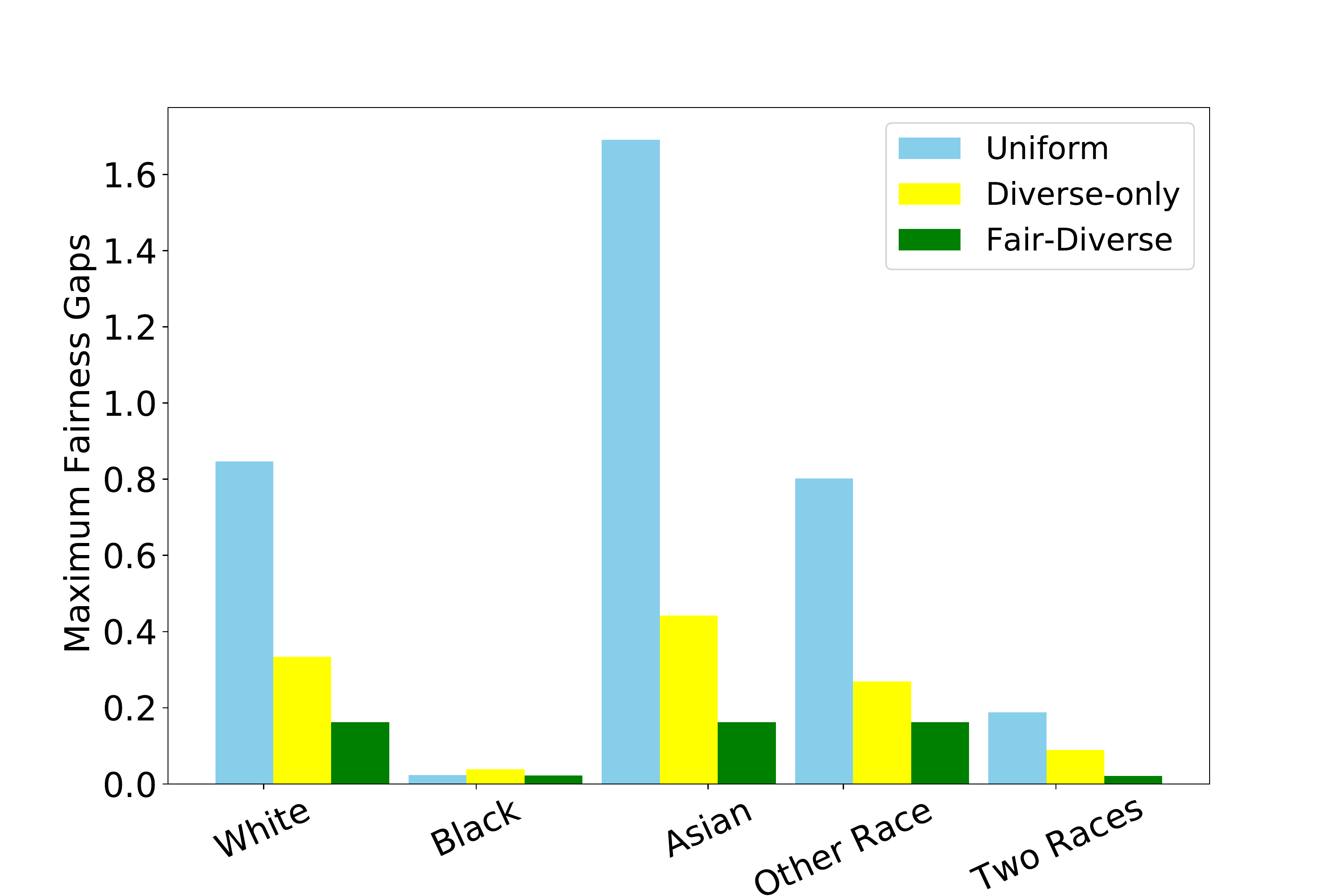}
  % \vspace{-5mm}
    \caption{\small Fairness gaps across different population subgroups}
    \label{fig:gaps2}
    \end{minipage}
    \begin{minipage}{0.35\textwidth}
        \includegraphics[width=\linewidth]{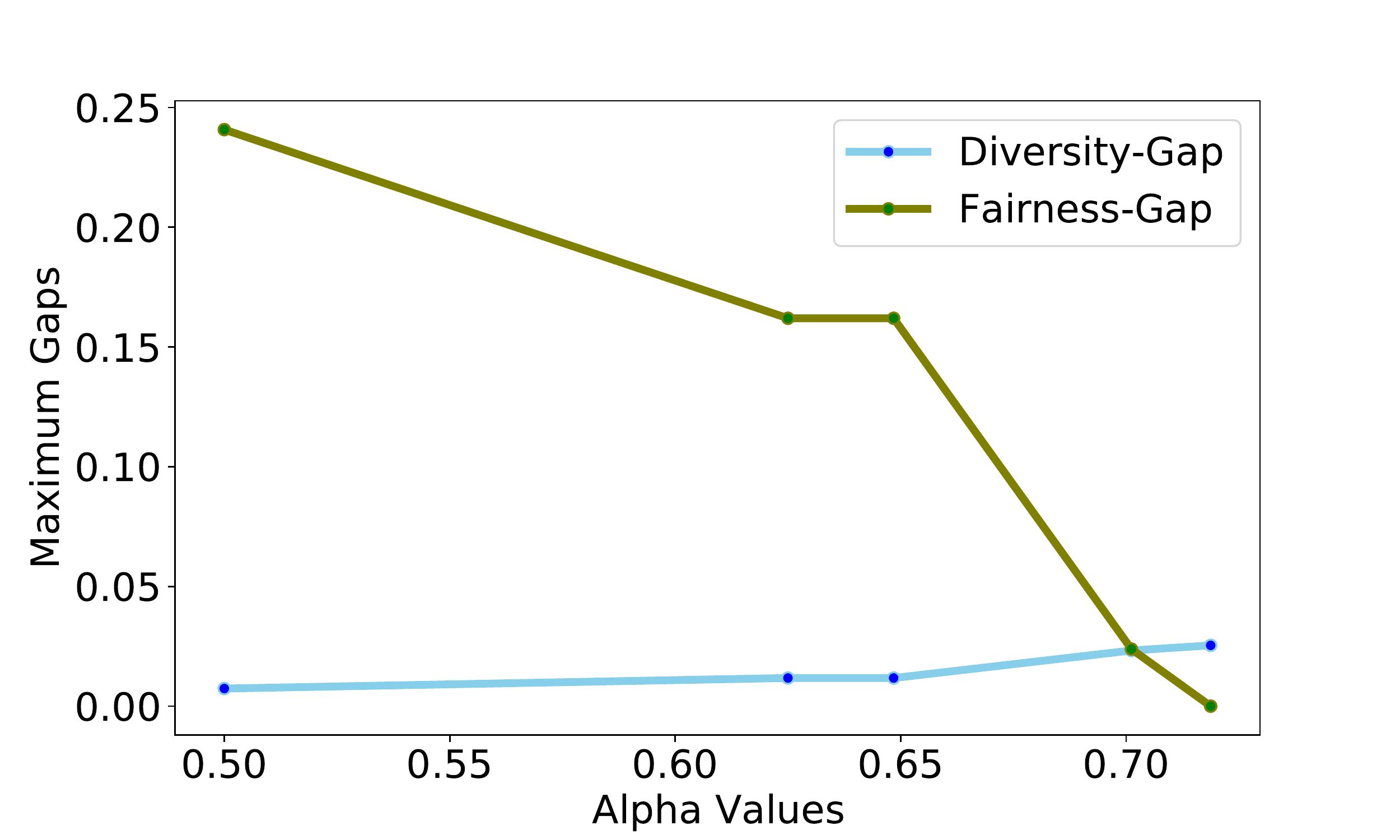}
    %\vspace{-5mm}
    \caption{\small Fairness and Diversity gaps based on different $\alpha$ values}
    \label{fig:gaps3}
    \end{minipage}
\end{figure}

Firstly, Figure~\ref{fig:gaps1} reveals that the \emph{Fair-Diverse} model reduces fairness and diversity gaps more compared to \emph{Diverse-only} and \emph{Fair-only} models. Although the \emph{Diverse-only} model eliminates the diversity gap, it fails to decrease the fairness gap. Similarly, the \emph{Fair-only} model eliminates the fairness gap but it fails to decrease the diversity gap.
Moreover, Figure~\ref{fig:gaps2} shows some considerable reduction in the gaps across different population subgroups (racial groups) using the \emph{Fair-Diverse} model. The comparison between the results with the \emph{Diverse-only} and the uniform allocation solution reveals the necessity of the \emph{Fair-Diverse} model in closing the gaps, and therefore, reducing the disparities across different population subgroups.

Lastly, Figure~\ref{fig:gaps3} illustrates the trade-off between fairness and diversity gaps considering different $\alpha$ values. Note that the $\alpha$ values represent the midpoint of the feasible range. We can observe that increasing the $\alpha$ value, which is the weight on the fairness component in $P2$, decreases the fairness gap as expected. However, the diversity gap will increase due to the fairness-diversity trade-off.  

Finally, as discussed in \S~\ref{sec:pof}, the Price of Fairness (PoF) can be evaluated using the difference of the fairness gap from the optimal solution of P2 that is obtained with and without any fairness constraints. POF could be defined as the fraction of the fairness gap that is obtained from the (\emph{Diverse-only}) allocation to the allocation solution based on the \emph{Fair-Diverse} model. If an allocation solution decreases the fairness gap more than the \emph{Fair-Diverse} model, the POF is less than one. Otherwise, PoF is > 1 and we will need to find a balance between the fairness and diversity objectives to decide on which allocation to choose. In the case of the Racial instance problem with $\alpha=0.71$, the PoF equals $\frac{0.4419}{0.025}=17.67$, which is significantly larger than 1.

%\nazanin{how to interpret}

%which means that we have reduced the gap by almost $\%$.

\section{Discussion}\label{sec:disc}

% \nazanin{The scope of application and limitation should be provided}
%\subsection{Estimating the Exposed Population of Intersectional Demographic Groups}\label{sec:intersectional} 

\new{In this section, we discuss the scope and limitation of the proposed fair-diverse resource allocation method, which is introduced in \S \ref{sec:technical}. 
Through the three case studies in \S \ref{sec:data}, we can see that the proposed method is applicable to the early stage of medical resource allocation when the resource is still considered to be scarce. 
As long as the required data are available and accurate, the proposed allocation scheme can be applied. 
We want to point out to some limitations of that the proposed approach can be further improved upon. 
}

\new{First, the exposure rate estimation can be improved. 
Recall that in \S \ref{sec:technical} the COVID-19 exposure rate for each intersectional group is denoted by $P(e|g_{i})$.
% To estimate these probabilities for intersectional population subgroups (e.g., <black,female,age30-39>), one can apply a Poisson regression model \cite{Poisson}, which is a widely used methodology when the response variable is a count. 
However, the primary challenge while estimating COVID-19 exposure rates for intersectional subgroups (e.g., <black,female,age30-39>), is the lack of intersectional population data in each Zipcode, denotes by $S_{ij}$ in \S \ref{sec:technical}. Given the availability of the required data format, one could use a Poisson regression model to obtain the intersectional exposure rates.}
% \textcolor{red}{From Lulu: this paragraph needs to be rewritten. Many duplications.}
\new{Second, the proposed method does not consider the profession of residents in prioritizing the sub-groups.
In the U.S., during the early vaccination process, priority was given to the sub-groups based on age, underlying health conditions, and professions. 
Front-line workers such as health care workers, grocery employees, K-12 educators, etc., whose works are essential to normal social functions were among the first one or two batches of vaccine receivers. 
In our work, we have estimated the COVID-19 exposure rate from the daily infected cases. 
It indirectly considers the residents' profession since front-line workers are expected to have a larger chance to be exposed to the virus.
But we are not able to directly include the profession of resident as a characteristic to define the sub-group of populations.
It is mainly due to the lack of data that has the job descriptions for residents (or the percentage) in each Zipcode area.
We plan to work with city officials or non-profit organizations that have more detailed population data and further improve our proposed method. 
}

Furthermore, estimating the intersectional population size could not be an appropriate implementation approach. Regarding the aforementioned data limitation issue, we solved the allocation problem for age, race, and gender as the sensitive attributes separately in \S\ref{sec:res1}. However, our model is able to return a fair-diverse allocation using the intersectional subgroups, given that the real intersectional population \new{is available for each area. In fact, applying the same framework on the intersectional data could be a future direction of this work.}

\new{Another future direction with respect to COVID-19 exposure rates, could incorporate a learning approach to obtain the rates. Several successful time-series techniques such as traditional ARIMA models or more novel approaches such as RNN-LSTM can be applied on the historical data. Time-series analysis would enable one to predict more precise exposure rates for the future time windows. }

Moreover, we only considered single treatment, mainly vaccination, however, the proposed approach can be generalized to incorporate multiple treatments. Our proposed fair-diverse allocation can be utilized at different stages of pandemic considering the updated exposure ratio. 

\new{In this paper, we mainly focused on the resource allocation to the centers and facilities. Although, designing allocation policy at a granular level of individuals from each center or facility is also critical, high-level policies have been shown to be more effective and feasible for deployment than individual levels ones. There are multiple barriers to the individual level deployment of the policies. For example, vaccine hesitancy \cite{Vaccinehesitancy-Covid1,Vaccinehesitancy-Covid2,Vaccinehesitancy-Covid3} is a widespread problem for many people and enforcing it could be in violation of social values such as Freedom of choice(FoC)\cite{FOC1,FOC2}. Moreover, designing strategies at individual level requires an extensive amount of data and medical considerations as vaccines(e.g. COVID19 vaccines) can cause serious allergic reaction and people with pre-existing allergies should avoid taking it\cite{Vaccine-allergy,Vaccine-allergy2,Vaccine-allergy3}.}
% ,...vaccine hesitations, ... \hadis{we need to find some articles discussing this}}

\new{Last but not least, applications of our proposed framework are beyond the allocation policy design use case. The solutions obtained for the allocation problem, could be helpful in identifying the vulnerable regions to incentive vaccination and for advertisement programs purposes. This will help maximizing the vaccination rate within higher risk communities in the future.}

\section{Conclusion}\label{sec:con}

In this paper, we propose the idea of fairness in scarce resource allocation problems (e.g, vaccine distribution) in terms of disparities across various population subgroups in different regions. To do so, we consider diversity and fairness components to design a Fair-Diverse allocation. We first formulate a general multi-objective (MO) problem $P1$, and propose the weighted sum method together with the LP relaxation to simplify it to $P2$. 
We then solve the LP-relaxed problem, $P2$, based on the fairness and diversity requirements as described in \S~\ref{sec:feasibility}. For this purpose, we propose a binary search approach, Algorithm \ref{alg:tune}, to find an optimal range for the trade-off parameter $\alpha$ in $P2$ such that the obtained solution is feasible.

Moreover, we have empirically analyzed our proposed methodologies in \S\ref{sec:exp} using COVID-19 datasets \new{ in three major and segregated US cities (New York City, Chicago, and Baltimore). We designed three instance problems for each city} based on different demographic attributes, race, age, and gender. We then implemented our fair-diverse model and compared it with other models (Diverse-only, Fair-only, and $\alpha=0.5$) to investigate the fairness criteria in the solution and highlight the necessity behind our approach. We have also discussed the fairness and diversity requirements, $\epsilon_{d}$ and $\epsilon_{f}$, and compared the Fair-Diverse allocation with other allocation solutions under different thresholds. Lastly, we have examined the price of fairness based on the associated gaps across different models. 

In brief, our empirical results reveal the paramount role of fairness criteria in decision-making problems involving scarce resource allocation (e.g, a vaccine allocation ). While certain minorities and population groups are more vulnerable to the COVID-19 virus, a Diverse-only (population-based) vaccine allocation can lead to higher fatality rates by neglecting the vulnerability of various population subgroups. We require to ensure a fair and diverse vaccine allocation to induce lower mortality rates across different regions. That is, we aim to find a \emph{decent balance} between the diversity and fairness measures in different geographical regions and allocate the resources accordingly.

\section{Acknowledgements}
L. Kang is partially supported by the National Science Foundation grant DMS-1916467.
\clearpage

\bibliographystyle{abbrv}
\bibliography{ref}

\begin{thebibliography}{10}

\bibitem{principal}
H.~Abdi and L.~J. Williams.
\newblock Principal component analysis.
\newblock {\em Wiley interdisciplinary reviews: computational statistics},
  2(4):433--459, 2010.

\bibitem{Covid-disparities}
D.~J. Alcendor.
\newblock Racial disparities-associated covid-19 mortality among minority
  populations in the us.
\newblock {\em Journal of clinical medicine}, 9(8):2442, 2020.

\bibitem{statisticatlas}
S.~Atlas.
\newblock The demographic statistical atlas of the united states.
\newblock
  \url{https://statisticalatlas.com/place/Illinois/Chicago/Race-and-Ethnicity},
  2018.
\newblock Accessed: 10-05-2020.

\bibitem{fairmlbook}
S.~Barocas, M.~Hardt, and A.~Narayanan.
\newblock Fairness and machine learning: Limitations and opportunities.
\newblock \url{fairmlbook.org}, 2019.

\bibitem{barocas2016big}
S.~Barocas and A.~D. Selbst.
\newblock Big data's disparate impact.
\newblock {\em Calif. L. Rev.}, 104:671, 2016.

\bibitem{fairness_old}
S.~K. Baruah, N.~K. Cohen, C.~G. Plaxton, and D.~A. Varvel.
\newblock Proportionate progress: A notion of fairness in resource allocation.
\newblock {\em Algorithmica}, 15(6):600--625, 1996.

\bibitem{becker1997optimal}
N.~G. Becker and D.~N. Starczak.
\newblock Optimal vaccination strategies for a community of households.
\newblock {\em Mathematical Biosciences}, 139(2):117--132, 1997.

\bibitem{Added-paper1}
H.~Beiki, S.~Mohammad~Seyedhosseini, V.~V.~Ponkratov, A.~Olegovna~Zekiy, and
  S.~A. Ivanov.
\newblock Addressing a sustainable supplier selection and order allocation
  problem by an integrated approach: a case of automobile manufacturing.
\newblock {\em Journal of Industrial and Production Engineering},
  38(4):239--253, 2021.

\bibitem{Covid-6cities}
J.~Benitez, C.~Courtemanche, and A.~Yelowitz.
\newblock Racial and ethnic disparities in covid-19: evidence from six large
  cities.
\newblock {\em Journal of Economics, Race, and Policy}, 3(4):243--261, 2020.

\bibitem{price}
D.~Bertsimas, V.~F. Farias, and N.~Trichakis.
\newblock The price of fairness.
\newblock {\em Operations research}, 59(1):17--31, 2011.

\bibitem{bertsimas2012efficiency}
D.~Bertsimas, V.~F. Farias, and N.~Trichakis.
\newblock On the efficiency-fairness trade-off.
\newblock {\em Management Science}, 58(12):2234--2250, 2012.

\bibitem{bertsimas2013fairness}
D.~Bertsimas, V.~F. Farias, and N.~Trichakis.
\newblock Fairness, efficiency, and flexibility in organ allocation for kidney
  transplantation.
\newblock {\em Operations Research}, 61(1):73--87, 2013.

\bibitem{bolton2003consumer}
L.~E. Bolton, L.~Warlop, and J.~W. Alba.
\newblock Consumer perceptions of price (un) fairness.
\newblock {\em Journal of consumer research}, 29(4):474--491, 2003.

\bibitem{managing}
J.~L. Bower.
\newblock {\em Managing the resource allocation process: A study of corporate
  planning and investment}.
\newblock Irwin Homewood, 1972.

\bibitem{chiu1989analysis}
D.-M. Chiu and R.~Jain.
\newblock Analysis of the increase and decrease algorithms for congestion
  avoidance in computer networks.
\newblock {\em Computer Networks and ISDN systems}, 17(1):1--14, 1989.

\bibitem{corbett2017algorithmic}
S.~Corbett-Davies, E.~Pierson, A.~Feller, S.~Goel, and A.~Huq.
\newblock Algorithmic decision making and the cost of fairness.
\newblock In {\em Proceedings of the 23rd acm sigkdd international conference
  on knowledge discovery and data mining}, pages 797--806, 2017.

\bibitem{whocandid}
A.~R. Cross.
\newblock Draft landscape of covid-19 candidate vaccines.
\newblock
  \url{https://www.redcrossblood.org/donate-blood/how-to-donate/how-blood-donations-help/blood-needs-blood-supply.html},
  2020.
\newblock Accessed: 04-15-2020.

\bibitem{demers1989analysis}
A.~Demers, S.~Keshav, and S.~Shenker.
\newblock Analysis and simulation of a fair queueing algorithm.
\newblock {\em ACM SIGCOMM Computer Communication Review}, 19(4):1--12, 1989.

\bibitem{Vaccine-allergy3}
A.~P. Desai, A.~P. Desai, and G.~J. Loomis.
\newblock Relationship between pre-existing allergies and anaphylactic
  reactions post mrna covid-19 vaccine administration.
\newblock {\em Vaccine}, 2021.

\bibitem{donahue2020fairness}
K.~Donahue and J.~Kleinberg.
\newblock Fairness and utilization in allocating resources with uncertain
  demand.
\newblock In {\em Proceedings of the 2020 Conference on Fairness,
  Accountability, and Transparency}, pages 658--668, 2020.

\bibitem{FOC2}
K.~Dowding and M.~Van~Hees.
\newblock Freedom of choice., 2009.

\bibitem{Vaccinehesitancy-Covid1}
A.~A. Dror, N.~Eisenbach, S.~Taiber, N.~G. Morozov, M.~Mizrachi, A.~Zigron,
  S.~Srouji, and E.~Sela.
\newblock Vaccine hesitancy: the next challenge in the fight against covid-19.
\newblock {\em European journal of epidemiology}, 35(8):775--779, 2020.

\bibitem{dubey2020representation}
R.~S. Dubey, G.~Laguzzi, and F.~Ruscitti.
\newblock On the representation and construction of equitable social welfare
  orders.
\newblock {\em Available at SSRN 3524071}, 2020.

\bibitem{dwork2012fairness}
C.~Dwork, M.~Hardt, T.~Pitassi, O.~Reingold, and R.~Zemel.
\newblock Fairness through awareness.
\newblock In {\em Proceedings of the 3rd innovations in theoretical computer
  science conference}, pages 214--226, 2012.

\bibitem{elzayn2019fair}
H.~Elzayn, S.~Jabbari, C.~Jung, M.~Kearns, S.~Neel, A.~Roth, and Z.~Schutzman.
\newblock Fair algorithms for learning in allocation problems.
\newblock In {\em Proceedings of the Conference on Fairness, Accountability,
  and Transparency}, pages 170--179, 2019.

\bibitem{emanuel2020fair}
E.~J. Emanuel, G.~Persad, R.~Upshur, B.~Thome, M.~Parker, A.~Glickman,
  C.~Zhang, C.~Boyle, M.~Smith, and J.~P. Phillips.
\newblock Fair allocation of scarce medical resources in the time of covid-19,
  2020.

\bibitem{feng2017multiobjective}
W.-H. Feng, Z.~Lou, N.~Kong, and H.~Wan.
\newblock A multiobjective stochastic genetic algorithm for the pareto-optimal
  prioritization scheme design of real-time healthcare resource allocation.
\newblock {\em Operations Research for Health Care}, 15:32--42, 2017.

\bibitem{emergency}
F.~Fiedrich, F.~Gehbauer, and U.~Rickers.
\newblock Optimized resource allocation for emergency response after earthquake
  disasters.
\newblock {\em Safety science}, 35(1-3):41--57, 2000.

\bibitem{Vaccine-allergy}
R.~E. Glover, R.~Urquhart, J.~Lukawska, and K.~G. Blumenthal.
\newblock Vaccinating against covid-19 in people who report allergies, 2021.

\bibitem{monopoly}
A.~C. Harberger.
\newblock Monopoly and resource allocation.
\newblock In {\em Essential readings in economics}, pages 77--90. Springer,
  1995.

\bibitem{hardt2016equality}
M.~Hardt, E.~Price, and N.~Srebro.
\newblock Equality of opportunity in supervised learning.
\newblock In {\em Advances in neural information processing systems}, pages
  3315--3323, 2016.

\bibitem{Kmeans}
J.~A. Hartigan and M.~A. Wong.
\newblock Algorithm as 136: A k-means clustering algorithm.
\newblock {\em Journal of the royal statistical society. series c (applied
  statistics)}, 28(1):100--108, 1979.

\bibitem{ho2019branch}
T.-Y. Ho, S.~Liu, and Z.~B. Zabinsky.
\newblock A branch and bound algorithm for dynamic resource allocation in
  population disease management.
\newblock {\em Operations Research Letters}, 47(6):579--586, 2019.

\bibitem{huaizhou2013fairness}
S.~Huaizhou, R.~V. Prasad, E.~Onur, and I.~Niemegeers.
\newblock Fairness in wireless networks: Issues, measures and challenges.
\newblock {\em IEEE Communications Surveys \& Tutorials}, 16(1):5--24, 2013.

\bibitem{jaffe1981bottleneck}
J.~Jaffe.
\newblock Bottleneck flow control.
\newblock {\em IEEE Transactions on Communications}, 29(7):954--962, 1981.

\bibitem{kelly1998rate}
F.~P. Kelly, A.~K. Maulloo, and D.~K. Tan.
\newblock Rate control for communication networks: shadow prices, proportional
  fairness and stability.
\newblock {\em Journal of the Operational Research society}, 49(3):237--252,
  1998.

\bibitem{koonin2020strategies}
L.~M. Koonin, S.~Pillai, E.~B. Kahn, D.~Moulia, and A.~Patel.
\newblock Strategies to inform allocation of stockpiled ventilators to
  healthcare facilities during a pandemic.
\newblock {\em Health security}, 18(2):69--74, 2020.

\bibitem{lan2010axiomatic}
T.~Lan, D.~Kao, M.~Chiang, and A.~Sabharwal.
\newblock {\em An axiomatic theory of fairness in network resource allocation}.
\newblock IEEE, 2010.

\bibitem{le2020covid}
T.~T. Le, Z.~Andreadakis, A.~Kumar, R.~G. Roman, S.~Tollefsen, M.~Saville, and
  S.~Mayhew.
\newblock The covid-19 vaccine development landscape.
\newblock {\em Nat Rev Drug Discov}, 19(5):305--306, 2020.

\bibitem{lum2016predict}
K.~Lum and W.~Isaac.
\newblock To predict and serve?
\newblock {\em Significance}, 13(5):14--19, 2016.

\bibitem{matrajt2013optimal}
L.~Matrajt, M.~E. Halloran, and I.~M. Longini~Jr.
\newblock Optimal vaccine allocation for the early mitigation of pandemic
  influenza.
\newblock {\em PLoS Comput Biol}, 9(3):e1002964, 2013.

\bibitem{Covid-disparities2}
L.~S. Mu{\~n}oz-Price, A.~B. Nattinger, F.~Rivera, R.~Hanson, C.~G. Gmehlin,
  A.~Perez, S.~Singh, B.~W. Buchan, N.~A. Ledeboer, and L.~E. Pezzin.
\newblock Racial disparities in incidence and outcomes among patients with
  covid-19.
\newblock {\em JAMA network open}, 3(9):e2021892--e2021892, 2020.

\bibitem{nesbitt2016increasing}
S.~Nesbitt and R.~E. Palomarez.
\newblock increasing awareness and education on health disparities for health
  care providers.
\newblock {\em Ethnicity \& disease}, 26(2):181, 2016.

\bibitem{chicagocity}
C.~of~Chicago.
\newblock City of chicago data portal.
\newblock
  \url{https://data.cityofchicago.org/browse?tags=covid-19&sortBy=relevance },
  2020.
\newblock Accessed: 09-15-2020.

\bibitem{Covid-Maryland}
M.~D. of~Health(MDH).
\newblock Coronavirus disease 2019 outbreak data.
\newblock \url{https://coronavirus.maryland.gov/ }, 2021.
\newblock Accessed: 07-07-2021.

\bibitem{ogryczak2007multicriteria}
W.~Ogryczak.
\newblock Multicriteria models for fair resource allocation.
\newblock {\em Control and Cybernetics}, 36(2):303--332, 2007.

\bibitem{ogryczak2005telecommunications}
W.~Ogryczak, M.~Pioro, and A.~Tomaszewski.
\newblock Telecommunications network design and max-min optimization problem.
\newblock {\em Journal of telecommunications and information technology}, pages
  43--56, 2005.

\bibitem{ong2021covid}
P.~M. Ong, C.~Pech, N.~R. Gutierrez, and V.~M. Mays.
\newblock Covid-19 medical vulnerability indicators: A predictive, local data
  model for equity in public health decision making.
\newblock {\em International Journal of Environmental Research and Public
  Health}, 18(9):4829, 2021.

\bibitem{CDPH}
C.~D. Portal.
\newblock Covid-19 daily cases, deaths, and hospitalizations, 2020.

\bibitem{Vaccine-allergy2}
A.~Remmel.
\newblock Covid vaccines and safety: what the research says, 2021.

\bibitem{Added-paper2}
C.~N. Rosyidi, S.~N. Hapsari, and W.~A. Jauhari.
\newblock An integrated optimization model of production plan in a large steel
  manufacturing company.
\newblock {\em Journal of Industrial and Production Engineering},
  38(3):186--196, 2021.

\bibitem{FOC1}
A.~Sen.
\newblock Freedom of choice: concept and content.
\newblock {\em European economic review}, 32(2-3):269--294, 1988.

\bibitem{singh2020fairness}
B.~Singh.
\newblock Fairness criteria for allocating scarce resources.
\newblock {\em Optimization Letters}, pages 1--9, 2020.

\bibitem{tanner2008finding}
M.~W. Tanner, L.~Sattenspiel, and L.~Ntaimo.
\newblock Finding optimal vaccination strategies under parameter uncertainty
  using stochastic programming.
\newblock {\em Mathematical biosciences}, 215(2):144--151, 2008.

\bibitem{tayfur2009model}
E.~Tayfur and K.~Taaffe.
\newblock A model for allocating resources during hospital evacuations.
\newblock {\em Computers \& Industrial Engineering}, 57(4):1313--1323, 2009.

\bibitem{Vaccinehesitancy-Covid3}
G.~Troiano and A.~Nardi.
\newblock Vaccine hesitancy in the era of covid-19.
\newblock {\em Public Health}, 2021.

\bibitem{Vaccinehesitancy-Covid2}
P.~Verger and E.~Dub{\'e}.
\newblock Restoring confidence in vaccines in the covid-19 era, 2020.

\bibitem{verweij2009moral}
M.~Verweij.
\newblock Moral principles for allocating scarce medical resources in an
  influenza pandemic.
\newblock {\em Journal of Bioethical Inquiry}, 6(2):159--169, 2009.

\bibitem{yao2017beyond}
S.~Yao and B.~Huang.
\newblock Beyond parity: Fairness objectives for collaborative filtering.
\newblock In {\em Advances in Neural Information Processing Systems}, pages
  2921--2930, 2017.

\bibitem{yarmand2014optimal}
H.~Yarmand, J.~S. Ivy, B.~Denton, and A.~L. Lloyd.
\newblock Optimal two-phase vaccine allocation to geographically different
  regions under uncertainty.
\newblock {\em European Journal of Operational Research}, 233(1):208--219,
  2014.

\bibitem{yip2021healthcare}
J.~Y.-C. Yip.
\newblock Healthcare resource allocation in the covid-19 pandemic: Ethical
  considerations from the perspective of distributive justice within public
  health.
\newblock {\em Public Health in Practice}, 2:100111, 2021.

\bibitem{zehlike2017fa}
M.~Zehlike, F.~Bonchi, C.~Castillo, S.~Hajian, M.~Megahed, and R.~Baeza-Yates.
\newblock Fa* ir: A fair top-k ranking algorithm.
\newblock In {\em Proceedings of the 2017 ACM on Conference on Information and
  Knowledge Management}, pages 1569--1578, 2017.

\bibitem{Covid-disparities3}
J.~Zelner, R.~Trangucci, R.~Naraharisetti, A.~Cao, R.~Malosh, K.~Broen,
  N.~Masters, and P.~Delamater.
\newblock Racial disparities in coronavirus disease 2019 (covid-19) mortality
  are driven by unequal infection risks.
\newblock {\em Clinical Infectious Diseases}, 72(5):e88--e95, 2021.

\bibitem{vzliobaite2017measuring}
I.~{\v{Z}}liobait{\.e}.
\newblock Measuring discrimination in algorithmic decision making.
\newblock {\em DATA MIN KNOWL DISC}, 31(4):1060--1089, 2017.

\end{thebibliography}
\vspace{200mm}

\appendix
\pagebreak
{\bf \LARGE APPENDIX}
\section{Proof of Proposition 1}
\noindent{\bf Proposition 2.}
{\it
% Let $X'\neq \emptyset$ be the set of feasible solution of $P3$. Let $A(\alpha)= (1-\alpha)\mathcal{D}+\alpha \mathcal{F})$. Given $x^* \in X'$, $\exists \alpha^* s.t. \{\argmin_x A(\alpha^*)|\sum_{j\in M} x_j\leq b, x_j \geq 0\} =x^*$.
Let be $X'=\{x|\mathcal{D}(x)<\epsilon_d, \mathcal{F}(x)<\epsilon_f, \sum_{j\in M} x_j=b, x_j \geq 0\}$:
\begin{enumerate}

    \item If $X'\neq \emptyset$: given $x^* \in X'$, $\exists \alpha^* $ such that $x^*$ is an optimal solution of $P2$.
    \item If $\not\exists ~\alpha^* $ such that $x^*$ is an optimal solution of $P2$ that satisfies both conditions $\mathcal{D}(x)<\epsilon_d$ and $ \mathcal{F}(x)<\epsilon_f$, then $X' = \emptyset$.
\end{enumerate}
% and $X'\neq \emptyset$. Given $x^* \in X'$, $\exists \alpha^* $ such that $x^*$ is an optimal solution of $P2$.
}
\begin{proof}
We first argue that the feasibility problem is equivalent to the following optimization model:
\begin{align*}
P':\\
& \min \mathcal{D}(x)\\
\nonumber &\mbox{ s.t. }\\
 &\mathcal{F}(x)= \epsilon\\
\quad &\mathcal{D}^+_j(x)\leq \mathcal{D}(x), \quad \forall j \in M\\
 & \mathcal{D}^-_j(x )\geq \mathcal{D}(x), \quad \forall j \in M\\
&\mathcal{F}^+_i(x)\geq \mathcal{F}(x), \quad \forall i \in I\\
& \mathcal{F}^-_i(x)\leq \mathcal{F}(x), \quad \forall i \in I\\
 &\sum_{j\in M} x_j = b\\
 &\;x_j \geq 0, \quad\quad \forall j \in M\\
% &\;x_j\in \mathbb{Z},\quad\quad \forall j \in M
\end{align*}
Consider $\epsilon\in [0,\epsilon_{\mathcal{F}}]$. Let $x^*$ be the optimal solution of $P'$. Hence, $\mathcal{F}(x^*)=\epsilon\leq \epsilon_{\mathcal{F}}$ satisfies the fairness constraint. Now, if $\mathcal{D}(x^*) \geq \epsilon_{\mathcal{D}}$, we conclude that for $\epsilon$ there does not exist a feasible solution that satisfies the diversity constraint. Otherwise, that feasible solution would minimize $P'$.

Introducing a Lagrangian multiplier $\lambda$, we now define the dual Lagrangian transformation of $P'$ to be:
\begin{align*}
P'':\\
& \min \mathcal{D}(x)+\lambda(\mathcal{F}(x)-\epsilon)\\
\nonumber &\mbox{ s.t. }\\
\quad &\mathcal{D}^+_j(x)\leq \mathcal{D}(x), \quad \forall j \in M\\
& \mathcal{D}^-_j(x )\geq \mathcal{D}(x), \quad \forall j \in M\\
&\mathcal{F}^+_i(x)\geq \mathcal{F}(x), \quad \forall i \in I\\
& \mathcal{F}^-_i(x)\leq \mathcal{F}(x), \quad \forall i \in I\\
&\sum_{j\in M} x_j= b\\
&\;x_j \geq 0, \quad\quad \forall j \in M\\
% &\;x_j\in \mathbb{Z},\quad\quad \forall j \in M
\end{align*}
Now, we argue that given $x^* \in X'$ there exists an $\alpha^*$ such that $x^*$ is an optimum of $P2$.
Let $\lambda^*$ be the multiplier when $x^*$ is an optimum of $P''$, then $\forall x''\neq x^*$:
\[
\mathcal{D}(x^*)+\lambda^*(\mathcal{F}(x^*)-\epsilon_{\mathcal{F}}) \leq \mathcal{D}(x'')+\lambda^*(\mathcal{F}(x'')-\epsilon_{\mathcal{F}})\]
\[
\implies \mathcal{D}(x^*)+\lambda^*\mathcal{F}(x^*)-\lambda^*\epsilon_{\mathcal{F}} \leq \mathcal{D}(x'')+\lambda^*\mathcal{F}(x'')-\lambda^*\epsilon_{\mathcal{F}}
\]
\[
\implies \mathcal{D}(x^*)+\lambda^*\mathcal{F}(x^*)\leq \mathcal{D}(x'')+\lambda^*\mathcal{F}(x'')
\]
Thus, $x^*$ is an optimum of $P2$ where $\lambda^*=\frac{(1-\alpha)}{\alpha}$.

On the other hand, if $\exists \alpha^*$ under which $x^*$ is optimal of $P2$, and $\mathcal{D}(x^*)<\epsilon_d$ and $\mathcal{F}(x^*)<\epsilon_f$, then $x^* \in X'$ and $X'\neq \emptyset$. This means if such $\alpha^*$ does not exist, there is no feasible solution under $\epsilon_d$ and $\epsilon_f$.
\end{proof}

\section{Proof of Proposition 2}

\noindent{\bf Proposition 2.}
{\it
Let $x_1^*$ and $x_2^*$ be optimum solutions of $P2$ given $\alpha_1$ and $\alpha_2$, respectively. 
It can be shown that if $\alpha_1\leq \alpha_2$, then $\mathcal{F}(x_2^*)\leq \mathcal{F}(x_1^*)$ and $\mathcal{D}(x_1^*)\leq \mathcal{D}(x_2^*)$. 
}
\begin{proof}
Let $\beta=\frac{\alpha}{(1-\alpha)}$. Given $x_1^*$ and $x_2^*$ corresponding to $\beta_1$ and $\beta_2$, where $\beta_2<\beta_1$:
\begin{equation*}
\mathcal{D}(x_1^*)+\beta_1\mathcal{F}(x_1^*)\leq\mathcal{D}(x_2^*)+\beta_1\mathcal{F}(x_2^*)
\end{equation*}
\begin{equation*}
\mathcal{D}(x_2^*)+\beta_2\mathcal{F}(x_2^*)\leq\mathcal{D}(x_1^*)+\beta_2\mathcal{F}(x_1^*)
\end{equation*}
Adding the above Equations, we will have:
\begin{equation*}
\beta_1\mathcal{F}(x_1^*)+\beta_2\mathcal{F}(x_2^*)\leq\beta_1\mathcal{F}(x_2^*)+\beta_2\mathcal{F}(x_1^*)
\end{equation*}
\[
\implies (\beta_2-\beta_1)(\mathcal{F}(x_2^*)-\mathcal{F}(x_1^*)) \leq 0
\]
which implies $\mathcal{F}(x_2^*)\leq\mathcal{F}(x_1^*)$. The monotonicity proof for $\mathcal{D}(x)$ is the same with $\beta=\frac{(1-\alpha)}{\alpha}$. If $\alpha_2 < \alpha_1$ then $\mathcal{D}(x_1^*)>\mathcal{D}(x_2^*)$.
\end{proof}

\new{\section{Exposure Rates}}

% \hadis{add a brief description for these tables}

\begin{figure}[h]  
\begin{minipage}[b]{0.47\textwidth}
    \centering
  \begin{adjustbox}{width=\columnwidth,center}
  \setlength\tabcolsep{1pt} 
    \begin{tabular}{ll}
    %\toprule
    \multicolumn{1}{p{13.25em}}{\textbf{Demographic Groups}} & \multicolumn{1}{p{7.915em}}{\textbf{Exposure rates}} \\
    \midrule
    Female & 0.095383 \\
    \midrule
    Male  & 0.09804 \\
    \midrule
    Age\_0\_4 & 0.035865 \\
    \midrule
    Age\_5\_12 & 0.057509 \\
    \midrule
    Age\_13\_17 & 0.061328 \\
    \midrule
    Age\_18\_24 & 0.081399 \\
    \midrule
    Age\_25\_34 & 0.106335 \\
    \midrule
    Age\_35\_44 & 0.107127 \\
    \midrule
    Age\_45\_54 & 0.106302 \\
    \midrule
    Age\_55\_64 & 0.126374 \\
    \midrule
    Age\_65\_74 & 0.130066 \\
    \midrule
    Age\_75+ & 0.173563 \\
    \midrule
    Asian/Pacific-Islander & 0.070976 \\
    \midrule
    Black/African-American & 0.063082 \\
    \midrule
    Hispanic/Latino & 0.207524 \\
    \midrule
    White & 0.049382 \\
    \end{tabular}
    \end{adjustbox}
    \caption{Exposure Rates-New York}
    \label{tab:rates_NYC}
    \end{minipage}
    \hfill
\begin{minipage}[b]{0.44\textwidth}
    \centering
  \begin{adjustbox}{width=\columnwidth,center}
  \setlength\tabcolsep{1pt} 
    \begin{tabular}{ll}
    %\toprule
    \multicolumn{1}{p{13.25em}}{\textbf{Demographic Groups}} & \multicolumn{1}{p{7.915em}}{\textbf{Exposure rates}} \\
       \midrule
    White & 0.050315 \\
    \midrule
    Black or African American & 0.083406 \\
    \midrule
    American Indian & 0.06383 \\
    \midrule
    Asian & 0.020939 \\
    \midrule
    Hawaiian or pacific Islander & 0.030405 \\
    \midrule
    Other race & 0.213587 \\
    \midrule
    Male  & 0.082152 \\
    \midrule
    Female & 0.085886 \\
    \midrule
    Age\_0\_9 & 0.042091 \\
    \midrule
    Age\_10\_19 & 0.056692 \\
    \midrule
    Age\_20\_29 & 0.092926 \\
    \midrule
    Age\_30\_39 & 0.120977 \\
    \midrule
    Age\_40\_49 & 0.087818 \\
    \midrule
    Age\_50\_59 & 0.086137 \\
    \midrule
    Age\_60\_69 & 0.099737 \\
    \midrule
    Age\_70\_79 & 0.087173 \\
    \midrule
    Age\_80+ & 0.07707 \\
    \end{tabular}
    \end{adjustbox}
  \caption{Exposure Rates-Baltimore}
  \label{tab:rates_BT}
\end{minipage}
\end{figure}

\end{document}